\documentclass[11pt]{article}
\usepackage{graphicx}
\usepackage{amsfonts}
\usepackage{amssymb}
\usepackage{amsmath,amsthm}
\usepackage{a4wide}
\usepackage{color}
\newcommand{\und}{\underline}
\newcommand{\dis}{\displaystyle}

 \newcommand{\nn}{\nonumber}

\numberwithin{equation}{section}
\theoremstyle{plain}
\newtheorem{theorem}{Theorem}[section]
\newtheorem{lemma}[theorem]{Lemma}
\newtheorem{corollary}[theorem]{Corollary}

\newtheorem{proposition}[theorem]{Proposition}
\theoremstyle{definition}
\newtheorem{definition}{Definition}[section]

\def\si{{\sigma}}

\def\ga{{\gamma}}
\def\om{{\omega}}
\def\Om{{\Omega}}

\def\N{{\mathbb N}}

\def\eps{{\epsilon}}

\def\square{\ifmmode\sqr\else{$\sqr$}\fi}
\def\sqr{\vcenter{
         \hrule height.1mm
         \hbox{\vrule width.1mm height2.2mm\kern2.18mm\vrule width.1mm}
         \hrule height.1mm}}                  

\newcommand{\col}[1]{\textcolor[rgb]{1,0,0}{#1}}

 \title{Hydrodynamic limit
in a particle system\\ with topological interactions}
\author{
Gioia Carinci$^{\textup{{\tiny(a)}}}$,
Anna De Masi $^{\textup{{\tiny(b)}}}$,\\
Cristian Giardin{\`a} $^{\textup{{\tiny(a)}}}$,
Errico Presutti $^{\textup{{\tiny(c)}}}$.\\\\
{\small $^{\textup{(a)}}$ Dipartimento di Dipartimento di Scienze fisiche, informatiche e matematiche,}\\
{\small Universit\`a di Modena e Reggio Emilia, via Campi 213/b, 41125 Modena, Italy}
\\
{\small $^{\textup{(b)}}$ Dipartimento di Ingegneria e Scienze dell'Informazione e Matematica,}\\
{\small Universit\`a di L'Aquila, via Vetoio 1, 67100 L'Aquila, Italy}\\
{\small $^{\textup{(c)}}$ GSSI, viale F. Crispi 7, 67100 L'Aquila, Italy}\\
{\small }\\
}

\date{\today}

\begin{document}

%
\maketitle

\begin{abstract}

We study a system of particles in the interval $[0,\eps^{-1}]\cap \mathbb Z$,
$\eps^{-1}$ a positive integer.  The particles move as
symmetric independent random walks
(with reflections at the endpoints);  simultaneously
new particles are injected at site $0$ at rate $j\eps$ ($j>0$)
and removed at same rate
from the rightmost occupied site.  The removal mechanism is therefore
of topological rather than metric nature.  The determination of the rightmost occupied site
requires a knowledge of the entire configuration and prevents
from using correlation functions techniques.

We prove using stochastic inequalities
that the system has a hydrodynamic limit, namely that
under suitable assumptions on the initial configurations, the law of the
density fields
$\eps \sum \phi(\eps x) \xi_{\eps^{-2}t}(x)$ ($\phi$ a test function, $\xi_t(x)$ the number of particles at  site  $x$ at time $t$)
concentrates in the limit $\eps\to 0$
on the deterministic value $\int \phi \rho_t$, $\rho_t$ interpreted as the limit density at time $t$.
We characterize the limit $\rho_t$ as a weak solution in terms of barriers
of a limit free boundary problem.

\end{abstract}

%
%
%
%
%
%
%
%
%
%
%



\vskip1cm

\section{Introduction and model definition}

This paper is inspired by the analysis in \cite{dfp} and we are
indebted to Pablo Ferrari for discussions and in particular
for suggesting the
inequalities in Section \ref{sesta}.  This is a first in a series of three papers
where we study a particle system whose hydrodynamic limit
is described by a free boundary problem.

Our system
is made  of particles confined to the lattice $[0,\eps^{-1}]\cap \mathbb Z$,
for brevity in the sequel we shall just write $[0,\eps^{-1}]$.  In this notation
$\eps^{-1}$ is a positive integer denoting the system size
and we will be eventually interested in the asymptotics as $\eps\to 0$.
The evolution is a Markov process $\{\xi_t,t\ge 0\}$ on the space of particles configurations
$\xi = (\xi(x))_{x\in [0, \eps^{-1}]}$, the component $\xi(x)\in\mathbb{N}$ is interpreted as
the number of particles at site $x$. The generator is denoted by
 \begin{equation}
\label{generatore}
L   = L^0   +  L_b  + L_{a}
\end{equation}
(the  dependence on $\eps$ is not made explicit).
$L^0$ is the generator
of the independent random walks process, it is
defined on  functions $f$ by
   \begin{equation}
\label{0.1}
L^0 f(\xi) =  \frac 12 \sum_{x=0}^{\eps^{-1}-1} L^0_{x,x+1} f(\xi)
    \end{equation}
    \begin{equation}
\label{generatore-bulk}
L^0_{x,x+1} f(\xi) = \xi(x) \left(f(\xi^{x,x+1}) -f(\xi)\right) + \xi(x+1) \left(f(\xi^{x+1,x}) -f(\xi)\right)
    \end{equation}
where
$\xi^{x,y}$ denotes the
configuration obtained from $\xi$ by removing
one particle from site $x$ and putting it at site $y$, i.e.
   \begin{equation*}
\xi^{x,y}(z) = \left\{
\begin{array}{ll}
\xi(z)  & \text{if } z \ne x,y,\\
\xi(z) -1  & \text{if } z = x,\\
\xi(z)+1  & \text{if } z=y.
\end{array} \right.
   \end{equation*}
Namely $L^0$ describes
independent symmetric random walks which jump
with equal probability
after an
exponential time of mean 1 to the nearest
neighbor sites, the jumps leading outside $[0,\eps^{-1}]$
being suppressed (reflecting boundary conditions).

The term $L_b$ in \eqref{generatore} is
   \begin{equation}
\label{generatore-birth}
L_{b} f(\xi) = j\eps \left(f(\xi^{+}) -f(\xi)\right),\quad \xi^{+}(x) = \xi(x) + \mathbf{1}_{x=0}
   \end{equation}
It describes the action of throwing into the system new particles at rate $\eps j$, $j>0$,
which then
land at  site  0; instead $L_a$ removes particles:
    \begin{equation}
\label{generatore-death}
L_{a}f(\xi) =  j\eps \left(f(\xi^{-}) -f(\xi)\right),\quad \xi^{-}(x) = \xi(x) - \mathbf{1}_{x=R(\xi)}
    \end{equation}
namely a particle is taken out from the edge  $R_\xi$ of the configuration $\xi$  defined as
   \begin{equation}
  \label{3.3.1.00}
 \text{ $R_\xi$ is such that:}\;\;\;  \begin{cases}
  \xi (y) >0 & \text{for $y= R_\xi$} \\
  \xi (y) =0 & \text{for $y>R_\xi$} \end{cases}
   \end{equation}
$L_{a}f(\xi)=0$ if $R_\xi$ does not exist, i.e.\ if $\xi\equiv 0$.

We interpret $L$ as the generator of a system
of independent walkers with ``current reservoirs'' which impose a positive
current $\eps j$ at  site  0 and at the edge
of the configuration.  See \cite{dptv1,dptv2} for a comparison with the
density reservoirs used in the analysis of the Fourier law.    Here
is a list of  the main issues which are studied in this and in the other papers in this series.

%
%
%
%

\medskip

\begin{itemize}

\item
The interaction described by {$L_a$} is highly non local as $R_\xi$ depends
on the  positions of all the particles.  This spoils any attempt to use
the BBGKY hierarchy of equations for the correlation
functions, as customary in $\eps$ perturbations of the independent system, see for instance \cite{dp1}.

\item
The $L_a$ interaction is ``topological rather than metric'',  {
as the influence on
a particle $i$ of a
particle $j$ only depends on whether $j$ is to the right or left of $i$ and not on
their distance}.
Topological interactions 
 appear  often in
natural sciences as in population dynamics, in particular
the motion of crowds of people {\cite{piccoli}},
or of animals {\cite{parisi}}.   Our result  shows that there are
natural examples
in physical systems as well.  The relative simplicity of our model allows
 a rigorous  analysis of such an interaction.

\item  To the left of $R_\xi$ the particles do not feel the $L_a$ interaction
and move freely, but $R_\xi$ depends  on the configuration of particles and hence on time  as well.
Ours therefore is a microscopic model for a free boundary problem and one
may thus guess  that the hydrodynamic limit
is also ruled by
a free boundary problem.  In such a case the
hydrodynamic equations would be the linear heat equation in an open, time dependent
space interval   with suitable
boundary conditions complemented by a law for the speed of the right boundary.

\item  The action of $L_b$ and $L_a$ is to add from
the left and respectively remove from the right  particles
at rate $\eps j$.  They act therefore as ``current reservoirs''
 \cite{dpt,dptv1,dptv2}  because they are imposing
a current $\eps j$  (recall that for density reservoirs  \cite{DF, CGGR}  the
particles current scales by $\eps$).  Supposing the validity of
Fick's law the stationary macroscopic profiles are then linear
functions with slope $-2j$: there are therefore infinitely many such profiles
(as here the boundary densities are not fixed).
Two scenarios are then possible:  either
there is a preferential profile or there is a second time scale beyond the hydrodynamical one,
where we see that such profiles are not stationary.

\end{itemize}

\medskip

\noindent
We shall give  answers  to most of the above issues, our
main results being stated in the next section.

\vskip2cm

\section{Main results}
\label{seconda}

Macroscopic profiles are functions $u \in L^\infty([0,1],\mathbb R_+)$ that
we also regard as positive Borel measures on $[0,1]$ via the correspondence
$u\to u\,dr$.  For any Borel positive measure $\mu$ on $[0,1]$ we define
	\begin{equation*}
F(r; \mu) = \int_r^1  \mu(dr'),\quad r\in [0,1] 
	\end{equation*}
setting, by an abuse of notation,
	\begin{equation}
	\label{ee8.1}
F(r; u) = \int_r^1  u(r')\,dr',\quad r\in [0,1]
	\end{equation}
We then say that   $u \in L^\infty([0,1],\mathbb R_+)$
has ``an edge'' $R(u)$ if
	\begin{equation}
	\label{ee8.1bis}
R(u)=\inf\{r: F(r; u) = 0\} <1
	\end{equation}
The definition extends naturally to Borel positive measures $\mu$ on $[0,1]$.

\medskip

\begin{definition} [Assumptions on the initial macroscopic profile]
   \label{defin0.0}

We denote by  $\rho_{\rm init}$  the initial macroscopic profile, we suppose
that   $\rho_{\rm init}\in L^\infty([0,1],\mathbb R_+)$.   \end{definition}

\medskip
\noindent {\bf Remark}. For some results we will
need  extra assumptions, namely that
$\rho_{\rm init} \in C([0,1],\mathbb R_+)$ and/or that  it has an ``edge''.

\medskip

We shall next discuss in which way particle systems and evolution of
macroscopic profiles are related.

\vskip.5cm

\subsection*{\em Hydrodynamic limit.}

%
%

\medskip

Particle  configurations  $\xi$
are elements of $\mathbb N^{[0,\eps^{-1}]}$
which may be regarded
as  positive measures $\mu_\xi$ on the real interval $[0,\eps^{-1}]$ by
setting
 \[
 \mu_\xi = \sum_{x=0}^{\eps^{-1}} \xi(x)D_x
 \]
where $D_x$, the Dirac delta at $x$, is the probability measure supported by
the point $x$.  Analogously to \eqref{ee8.1} we set
	\begin{equation}
	\label{ee8.1.0}
F_\eps(x;\xi) = \int_x^{\eps^{-1}} \mu_\xi(dx')= \sum_{y\ge x} \xi(y),\quad x\in [0,\eps^{-1}]
	\end{equation}
and, as for the macroscopic profiles, we say that $\xi$ has an edge $R_\xi$ if
	\begin{equation}
	\label{ee8.1bis.0}
R_\xi=\inf\{x: F(x; \xi) = 0\} <\eps^{-1}
	\end{equation}
which means that $R_\xi<\eps^{-1}$ is the largest integer $x$ such that $\xi(x)>0$, in agreement
with \eqref{3.3.1.00}.
To compare macroscopic profiles and particles configurations we shall use the functions
$F_\eps(x;\xi)$ and $F(r; u)$. We define in particular the  local  averages:
 \begin{equation}
  \label{3.3.0}
 \mathcal A_\ell(x,\xi) := \frac 1{\ell} \Big( F_\eps(x;\xi)- F_\eps(x+\ell-1; \xi)\Big)=
  \frac 1{\ell} \sum_{y=x}^{ x+\ell-1} \xi(y)
\end{equation}
with  $\ell$ a positive integer and  $x\in [0,\eps^{-1}-\ell+1]$.  The corresponding
quantity for macroscopic profiles $u\in L^\infty([0,1],\mathbb R_+)$ is
 \begin{equation}
  \label{3.3.1}
  \mathcal A'_\ell(x,u)= \frac 1{\eps\ell} \Big( F(\eps x; u)-  F(\eps (x+\ell); u) \Big)
\end{equation}

%
%
%
%
%

\medskip

\begin{definition} [Assumptions on the initial particle configuration]
   \label{defin0.1}
We fix  $b<1$  suitably close to 1 and $a>0$ suitably small,
for the sake of definiteness we set $b=9/10$
and $a=1/20$.
We then denote by $\ell$ the integer part of $\eps^{-b}$
and suppose that for any $\eps$ the initial configuration $\xi$ verifies
 \begin{equation}
  \label{3.3.1}
\max_{x\in [0,\eps^{-1}-\ell+1]}\Big| \mathcal A_\ell(x,\xi)-
  \mathcal A'_\ell(x,\rho_{\rm init})\Big| \le \eps^a
\end{equation}
and moreover that if $\rho_{\rm init}$ has an edge $R(\rho_{\rm init})$, see   \eqref{ee8.1bis},
then
 \begin{equation}
  \label{3.3.1.00.bis}
  |\eps R_\xi-R(\rho_{\rm init})|\le \eps^a
   \end{equation}
with $R_\xi$ as in \eqref{ee8.1bis.0}.
We shall denote by $P^{(\eps)}_\xi$ the law of the  process $\{\xi_t,t\ge 0\}$
in the interval $[0,\eps^{-1}]$ with generator $L$ given in \eqref{generatore}
and started at time 0 from a configuration  $\xi$ as above.

\end{definition}

\medskip

Thus the initial configuration $\xi$ converges to $\rho_{\rm init}$ as $\eps\to 0$ in the sense of \eqref{3.3.1}.
Our first result proves that the convergence extends to all positive times  (but
in a weaker sense).

\medskip

  \begin{theorem}[Existence of hydrodynamic limit]
   \label{Teo:Hydro}

Let $\rho_{\rm init}\in L^\infty([0,1],\mathbb R_+)$ 
and  $\xi$
as in Definition \ref{defin0.1}.
Then there exists a non negative, continuous function $\rho(r,t)$,
$t> 0$, $r\in [0,1]$, such that
for any $r\in [0,1]$
   \begin{equation}
   \label{eq:Hydrobis}
\lim_{t \to 0} F(r;\rho(\cdot,t))  = F(r;\rho_{\rm init}(\cdot))
    \end{equation}
and such that for any $t>0$ and $\zeta>0$
   \begin{equation}
   \label{eq:Hydro}
\lim_{\eps \to 0}P^{(\eps)}_\xi\Big[\max_{x\in[0,\eps^{-1}]} |\eps F_\eps(x;\xi_{\eps^{-2}t}) - F(\eps x;\rho(\cdot,t))| \le \zeta \Big] = 1 
    \end{equation}
Moreover, if $\rho_{\rm init}\in C([0,1],\mathbb R_+)$ then $\rho(r,t)$ is continuous in $[0,1]\times \{t\ge 0\}$
and $\rho(r,0)=\rho_{\rm init}$.

   \end{theorem}

   \vskip.4cm

The above convergence implies weak convergence of the density fields against smooth test functions  $\phi$:
    $$
\lim_{\eps\to 0} P^{(\eps)}_\xi\left[\Big|\eps \sum_{x} \xi_{\eps^{-2}t}(x) \phi( \eps x) - \int_{0}^1 \phi(r)\rho(r,t) dr\Big| \le \zeta\right] = 1,\;\;\text{for all $\zeta>0$.}
   $$

\vskip.5cm

\subsection*{\em  The free boundary problem.}

Theorem \ref{Teo:Hydro} states  the  existence
and some regularity
properties of the hydrodynamic limit,
but does not say about its qualitative features: in particular
which equation is satisfied by the limit and which equation rules the motion of the
edge, if it exists.  The continuum analogue of our particle evolution is
 \begin{equation}
\label{macro.2}
\frac{\partial \rho}{\partial t}   = \frac 12
\frac{\partial^2 \rho}{\partial r^2} + j D_0 - j D_{R_t},\quad j>0
\end{equation}
where the first term (on the right hand side) corresponds to the random walk evolution,
$jD_0$ to the addition of particles at the origin and $j D_{R_t}$ to  the
removal of the rightmost particles.

In \cite{CDGP2} a suitable notion of quasi-solutions for
\eqref{macro.2} in $\mathbb R_+$ is given and it is proved that the limit of such  quasi-solutions  coincide
with the hydrodynamic limits found in Theorem \ref{Teo:Hydro}.  The main ingredient
in the proof is established here and it is based on the notion of upper and lower
barriers.   These are
``approximate  solutions''  of \eqref{macro.2} which bound from  below and from above  the hydrodynamic limit
$\rho(r,t)$, the inequalities being in the sense of mass transport.

This is defined as follows: two positive Borel measures $\mu$ and $\nu$ on $[0,1]$ are ordered
 with $\mu \le \nu$ if
    \begin{equation*}
F(r; \mu)    \le F(r; \nu) \;\;\;\text{for all $r\in [0,1]$\;.}
   \end{equation*}
 We shall apply the notion to measures in $\mathcal U$ defined as follows:

\medskip

   \begin{definition}
\label{defin:ee8.1}
 (The set $\mathcal U$ and the partial order).   $\mathcal U$ is the set of all positive Borel measures $u$ on $[0,1]$
 which have the form $u=c_uD_0 + \rho_u(r)dr$, $c_u\ge 0$, $\rho_u\in L^\infty([0,1],\mathbb R_+)$.  By an abuse
 of notation we shall also write the elements of $\mathcal U$ as $u=c_uD_0 +\rho_u$.
For any $u,v \in \mathcal U$
we then set
   \begin{equation}
\label{e4.3.0.2.1.000}
u\le v \;\;\;\text{iff}\;\;\;
F(r; u)    \le F(r; v) \;\;\;\text{for all $r\in [0,1]$\;.}
   \end{equation}
We also write $|u-v| = |c_u-c_v|   D_0 +|\rho_u-\rho_v|\in \mathcal U$ so that
      \begin{equation}
\label{e4.3.0.2.1.000.1}
|u-v| _1:=F(0;| u-v|)  = |c_u-c_v| + \int_0^1 |\rho_u(r)-\rho_v(r)|\, dr
   \end{equation}
is the total variation of the measure $u-v$.

\end{definition}

\medskip

 %
%
%


\medskip

\medskip
\begin{definition}
 (The cut and paste operator).
We define for any $\delta>0$ the subset
 $\mathcal U_\delta \subset \mathcal U$ as 
  \begin{equation}
\label{e4.3.0.0}
\mathcal U_\delta:=\{u= c_u D_0+\rho_u:\; F(0;\rho_u) > j\delta\}
   \end{equation}
and the {\em cut-and-paste} operator
$K^{(\delta)}:  \mathcal U_\delta\to \mathcal U$
   \begin{equation}
\label{e4.3.0.1}
K^{(\delta)} u = j\delta D_0+
\mathbf 1_{r\in [0,R_\delta(u)]}u,\quad R_\delta(u)=\inf\{r: F(r;u)=j\delta\}
   \end{equation}
   Observe that $F(0; K^{(\delta)} u )=F(0;u)$.

In the following definition of barriers we use
the Green function $G_\delta^{\rm neum}(r,r')$
(for the heat equation in $[0,1]$ with Neumann boundary conditions):
\begin{equation}
\label{e4.3}
G_t^{\rm neum}(r,r') = \sum_k G_t(r,r'_k),\quad
 G_t(r,r') = \frac{e^{-(r-r')^2/2t}}{\sqrt{2\pi t}}
\end{equation}
$r'_k$ being the images of $r'$ under
repeated reflections of the interval $[0,1]$ to its right and left
(see for instance \cite{KS} pag. 97 for details).
\end{definition}

We denote by $$G_\delta^{\rm neum} * f(r) = \int G_\delta^{\rm neum}(r,r')f(r')\,dr'$$
and observe that $F(0; G_\delta^{\rm neum} *u )=F(0;u)$ and
$G_\delta^{\rm neum} *u \in L^\infty([0,1];\mathbb R_+)$.
\medskip

\begin{definition} [Barriers]
\label{defin:e2.3}
Let $u\in L^\infty([0,1],\mathbb R_+)$ be such that $F(0; u) >0$. Then for all $\delta$ small enough   $u\in \mathcal U_\delta$
and for such $\delta$ we define
the ``barriers''  $S_{n\delta}^{(\delta,\pm)}(u) \in {\mathcal U}_{\delta}$,  $n\in \mathbb N$,
as follows: we set
$S^{(\delta,\pm)}_{0}(u)=u$, and, for $n\ge 1$,
\begin{eqnarray}
\label{e2.10}
&&
S^{(\delta,-)}_{n\delta}(u)= K^{(\delta)} G_\delta^{\rm neum} *S^{(\delta,-)}_{(n-1)\delta}(u)
\\&&
S^{(\delta,+)}_{n\delta}(u)= G_\delta^{\rm neum} * K^{(\delta)} S^{(\delta,{+})}_{(n-1)\delta}(u)
\nn
\end{eqnarray}
The families
$\{S_{n\delta}^{(\delta,+)}(u)\}_{\delta>0}$
are called upper barriers and $\{S_{n\delta}^{(\delta,-)}(u)\}_{\delta>0}$ lower barriers.

\end{definition}

\medskip
The functions $S_{n\delta}^{(\delta,\pm)}$ are obtained by alternating
the map $G^{\rm neum}_{\delta}$ (i.e.\ a diffusion)
and the cut and paste map $K^{(\delta)}$
(which takes out a mass $j\delta$ from the right and put it back
at the origin,  the macroscopic counterpart of $L_b$ and $L_a$).
It can be easily seen that unlike the original process $\xi_t$ the evolutions
$S_{n\delta}^{(\delta,\pm)}$ conserve the total mass,
that $S_{n\delta}^{(\delta,+)}$ maps $L^\infty$ into $C^\infty$ while  $S_{n\delta}^{(\delta,-)}$
has a singular component ($j\delta D_0$)
plus a $L^\infty$ component (which is $C^\infty$ inside its support).

\medskip

The name ``upper and lower barriers'' is justified by the following theorem:
%

\medskip

           \begin{theorem}[Separated classes]
   \label{prop-class}
Let $u\in \L^\infty([0,1],\mathbb R_+)$, $F(0;u)>0$, then 
\begin{equation}
\label{e4.3.0.2.1.00}
S_{t}^{(\delta,-)}(u) \le S_{t}^{(\delta',+)}(u) \quad
\text{for all $\delta,\delta',t$ such that $u \in \mathcal U_\delta \cap \mathcal U_{\delta'}$ and $t=k\delta=k'\delta'$, with $k,k'\in\N$}
\end{equation}
where the inequality is in the sense of Definition \ref{defin:ee8.1}.
          \end{theorem}

\medskip

It thus looks natural to look for elements which separate the barriers:

\medskip

\begin{definition}[Separating elements]
\label{defin:e2.3.00}
For a given non negative $u \in L^\infty$, the function $u=u(r,t)$, $r\in [0,1]$, $t\ge 0$, is below the upper barriers $\{S_{n\delta}^{(\delta,+)}(u)\}$
if
\begin{equation}
\label{e4.3.0.2.1.00.1}
 u(\cdot,t) \le S_{t}^{(\delta,+)}(u)(\cdot) \quad
\text{for all $\delta>0$ and $t$ such that $t=k\delta$,  $k\in \mathbb N$}
\end{equation}
It is above
the lower barriers $\{S_{n\delta}^{(\delta,-)}(u)\}$
if
\begin{equation}
\label{e4.3.0.2.1.00.1bis}
 u(\cdot,t) \ge S_{t}^{(\delta,-)}(u)(\cdot) \quad
\text{for all $\delta>0$ and $t$ such that $t=k\delta$,  $k\in \mathbb N$}
\end{equation}
If it is both above $\{S_{n\delta}^{(\delta,-)}(u)\}$ and below $\{S_{n\delta}^{(\delta,+)}(u)\}$
then $u(\cdot,t)$ separates the  barriers $\{S_{n\delta}^{(\delta,\pm)}(u)(\cdot)\}$.

%
Observe that if $u(\cdot,t)$   separates $\{S_{n\delta}^{(\delta,\pm)}(u)\}$ then $u(\cdot,0) = u(\cdot)$.

\end{definition}

\vskip.5cm

   \begin{theorem}
   [Existence and uniqueness of separating elements]
   \label{Teo:EU}
Let $u\in L^\infty([0,1],\mathbb R_+)$ and $F(0;u)>0$.  Then there exists
a unique function $u(r,t)$ which separates the barriers
$\{S_{n\delta}^{(\delta,\pm)}(u)\}$. $u(r,t)$ is continuous on the compacts
of $[0,1]\times (0,\infty)$ and $u(\cdot,t)$ converges weakly to $u(\cdot)$
as $t\to 0$.

   \end{theorem}

\medskip
More properties of the separating elements are established in Section \ref{ottava}, in particular we show that
they can be obtained as monotonic limits of the upper or the lower barriers.

\vskip.5cm

   \begin{theorem}[Characterization of hydrodynamic limit]
   \label{Teo:Char}
The hydrodynamic limit $\rho(r,t)$ of Theorem \ref{Teo:Hydro} separates
the barriers   $\{S_{n\delta}^{(\delta,\pm)}(\rho_{\rm init})\}$.
 %

   \end{theorem}

\vskip.5cm

\subsection*{\em  Super-hydrodynamic limit and further results.}

In {\cite{CDGP3} we shall study the stationary solutions of \eqref{macro.2}, they
are linear functions with slope $-2j$.  We shall prove that any weak solution (in the sense of barriers)
converges as $t\to \infty$ to a linear profile, the one with the same total mass as the initial state.
We shall also prove that at super-hydrodynamic times, i.e.\ times of order $\eps^{-3}$ the particle processes
is ``close'' to the manifold of linear profiles performing a brownian motion on such a set.

We conclude the list of results in this paper
by  a last theorem where we  identify the  limit equation for $\rho(\cdot,t)$ when $\rho_{\rm init}(\cdot)$ has no edge:

\medskip

  \begin{theorem}[{Hydrodynamic limit in the absence of an edge}]
   \label{TEO:2}
Let $\rho_{\rm init}$ such that $F(r; \rho_{\rm init}) \ge \alpha(1-r)$, $\alpha>0$, then there exists $T>0$ such that
$\rho(1,t)>0$ for $t\in [0,T]$ and $\rho(r,t)$ is given by
   \begin{equation}
   \label{e2.5}
 \rho(r,t) = G^{\rm neum}_t * \rho_{\rm init} (r) + j \int_0^t \{G^{\rm neum}_s (r,0) - G^{\rm neum}_s (r,1)\}ds,\;\; t\in [0,T]
    \end{equation}
$G^{\rm neum}_t(r,r')$ being the Green function of the heat equation in $[0,1]$ with Neumann conditions, see
\eqref{e4.3}.

   \end{theorem}

\vskip.5cm

\subsection*{\em  Strategy of proof.}

The key observation is that
if we anticipate/posticipate the addition and
removal of the particles which occur
in the true process in a given time interval then
we stochastically increase/decrease the final configuration
(in the sense of mass transport to the right, i.e. the microscopic version of \eqref{e4.3.0.2.1.000}).

To implement this we introduce the processes $\xi^{(\delta,\pm)}_{k\eps^{-2}\delta}, k \in \mathbb N$.
If for the true process
the number of added and removed particles in the time interval
$[k\eps^{-2}\delta,(k+1)\eps^{-2}\delta]$
is equal to $N_{k;\pm}$ then $\xi^{(\delta,-)}_{(k+1)\eps^{-2}\delta}$ is obtained from
$\xi^{(\delta,-)}_{k\eps^{-2}\delta}$ by letting it evolve with generator $L^0$ and at the end
adding $N_{k;+}$ particles at 0 and then  removing the rightmost $N_{k;-}$ particles.  In a similar fashion $\xi^{(\delta,+)}_{(k+1)\eps^{-2}\delta}$ is obtained by reversing the order of the operations: first
the addition/removal and then after the free evolution.  We then have for all
$\delta>0$ and all  $k\in \mathbb N$
   \begin{equation}
  \label{2.14bb}
  \xi^{(\delta,-)}_{k\eps^{-2}\delta} \le \xi_{k\eps^{-2}\delta} \le \xi^{(\delta,+)}_{k\eps^{-2}\delta}\;\;\;\text{\rm stochastically}
     \end{equation}
(see Section \ref{sesta} for details,  in particular the definition of microscopic notion
of partial order).

The probabilistic part of the paper is essentially concentrated in the analysis of the
hydrodynamic limit of the process $\xi^{(\delta,\pm)}_{k\eps^{-2}\delta}$: in Section \ref{quarta}
we prove that it converges to $S^{(\delta,\pm)}_{k\delta}(u)$ (if the initial $\xi$ ``approximates''
$u$) where convergence is in the sense of \eqref{eq:Hydro}.  This is important because
it implies that the inequalities are preserved in the limit.

The hydrodynamic limit for the independent random walks process is
easy and well known in the literature, but in our case there is
an extra difficulty related to a macroscopic occupation at the origin, $\xi(0)\approx \eps^{-1}$,
due to the cut and paste operations.  This severely limits the
choice of the parameters ($b$ close to 1,
$a$ close to 0  which in normal situations have a much larger range of values) but luckily
some room is left.  Instead the convergence of the microscopic cut and paste to its
macroscopic counterpart is easy, as the variables $N_{k;\pm}$ are modulo negligible
deviations independent Poisson variables with mean $j\eps^{-1}\delta$.

Once we have convergence to  $S^{(\delta,\pm)}_{k\delta}(u)$  we are left with the
analytic problem of studying the limits of the latter as $\delta\to 0$.  We first prove some regularity
properties uniform in $\delta$, see Section \ref{settima}, and then complete the proof of all theorems.

\vskip.5cm

\subsection*{\em  Sections content.}

In Section \ref{terza} we introduce the $\delta$-approximate
processes $\{\xi^{(\delta,\pm)}_t\}$ and prove that the law of the total particles number process
$|\xi_t|$ is a symmetric
random walk on $\mathbb N$ with reflection at the origin (a result which follows directly
from the definition of the process $\xi_t$).  We then state some consequences of such a result which will be used
in the sequel.

In Section \ref{quarta} we prove that if the initial configuration $\xi$
approximates a profile $u\in \mathcal U$ then
$\xi_{\eps^{-2}k\delta}^{(\delta,\pm)}$  converges
in law to
$S_{k\delta}^{(\delta,\pm)}(u)$ as $\eps\to 0$.
The proof
exploits duality for the independent process but
is not a consequence of well known results on
the hydrodynamic limit for independent particles
because we need to take into account the case when
there is a macroscopic occupation number at the origin.
As a consequence the bounds are not as strong as those
which appear in the literature.

In Section \ref{quinta} we introduce a probability space $(\Om,P)$ where we
can realize simultaneously all the
processes $\xi_t$ and $\xi_{\eps^{-2}k\delta}^{(\delta,\pm)}$ for all $\eps$.

In Section \ref{sesta} we  relate the true process $\xi_{\eps^{-2}k\delta}$ and the auxiliary ones $\xi_{\eps^{-2}k\delta}^{(\delta,\pm)}$
by stochastic
inequalities, in the sense of mass transport theory, exploiting
the realization of the process
of Section \ref{quinta}.  By using the convergence proved in Section \ref{quarta}
the inequalities extend to flows $S_{k\delta}^{(\delta,\pm)}$, thus proving
Theorem \ref{prop-class}.

In Section \ref{settima} we prove regularity properties of the
flows $S_{k\delta}^{(\delta,\pm)}$ which are uniform in $\delta$.

In Section \ref{ottava} we prove we first prove existence and
uniqueness of the separating element of barriers (Theorem \ref{Teo:EU})
and then deduce our main results
(Theorems  \ref{Teo:Hydro}
and \ref{Teo:Char}). We conclude by giving the
proof of  Theorem \ref{TEO:2}.

\vskip1cm

\section{The $\delta$-approximate  particle processes}
\label{terza}

In this Section we define the  stochastic processes $\xi^{(\delta,\pm)}_{k\eps^{-2}\delta}, k \in \mathbb N$ which are analogous to the barriers $S^{(\delta,\pm)}_{k\delta}$ of Definition \ref{defin:e2.3}. As we shall explain below, these processes are defined in such a way that the number of added and removed particles in the time interval
$[k\eps^{-2}\delta,(k+1)\eps^{-2}\delta]$, denoted by $N_{k;\pm}$, are the same as those in  the true process $\{\xi_t\}$.

The variables $N_{k;\pm}$, $k \in \mathbb N$ are determined by
the increments of process  $|\xi_t|$ yielding the particles' number
at time $t$. This last process,
despite the complexity of the full process $\xi_t$,
 is very simple:

\vskip.5cm

\begin{theorem}[Distribution of the particles' number]
\label{thme3.6a}
$|\xi_t|$ has the law
of a random walk {$(n_t)_{t\ge 0}$} on
$\mathbb N$ which jumps with equal probability by $\pm 1$
after an  {exponential} time of  {parameter} $2j\eps$,
the jumps leading to $-1$ being suppressed.

\end{theorem}

\medskip

\noindent
{\bf Proof.}  For any bounded function $f$ on $\mathbb N$ we have
\begin{equation}\label{Gen}
Lf(|\xi|) = j\eps \Big \{\big( f(|\xi|+1)-f(|\xi|)\big) + \mathbf 1_{|\xi|>0}
\big( f(|\xi|-1)-f(|\xi|)\big)\Big \}
\end{equation}
which  coincides with  the action of the generator of the random walk  {$(n_t)_{t\ge 0}$} on
the function $f(n)$.  This proves that the law of $|\xi_t|$ is the same as
that of the random walk.
\qed

\medskip

To introduce the $\delta-$approximate process we define
	\begin{equation}
	\label{e3.1aa}
N_{k;+} = \;\;\text{number of {upwards} jumps
of $|\xi_t|$ for $t\in [k\eps^{-2}\delta,(k+1)\eps^{-2}\delta]$}
	\end{equation}
	\begin{equation}
	\label{e3.2aa}
N_{k;-} =\;\;\text{number of  {downwards} jumps
of $|\xi_t|$ for $t\in [k\eps^{-2}\delta,(k+1)\eps^{-2}\delta]$}
	\end{equation}

	\medskip

\begin{definition}[{The $\delta-$approximated processes}]
\label{define3.1aa}
The processes $\xi_t^{(\delta,\pm)}$ are defined successively in
the time intervals $[k\eps^{-2}\delta,(k+1)\eps^{-2}\delta]$, $k\ge 0$.  We first distribute
the variables $N_{k;\pm} $
as the increments of the Markov process $(|\xi_t|)_{t\ge 0}$ starting from $|\xi_0^{(\delta,\pm)}|$.
Given such variables we use an  induction procedure  and suppose  $\xi_{k\eps^{-2}\delta}^{(\delta,-)} = \xi$
given.  Then $\xi_t^{(\delta,-)}$, $t \in [k\eps^{-2}\delta,(k+1)\eps^{-2}\delta)$ has
the law of the process $\xi^0_t$ with generator $L_0$  defined in \eqref{generatore-bulk}  starting from $\xi$ at time $k\eps^{-2}\delta$.
$\xi_{(k+1)\eps^{-2}\delta}^{(\delta,-)}$
is then
obtained from $\xi_{(k+1)\eps^{-2}\delta}^{0}$ by adding $N_{k;+} $ particles all at the origin
and then removing the $N_{k;-} $ rightmost particles.

$\xi_t^{(\delta,+)}$, $t\in (k\eps^{-2}\delta,(k+1)\eps^{-2}\delta]]$, is defined as the independent  random walk evolution
starting at time $k\eps^{-2}\delta$ from $\xi'$: $\xi'$ is obtained from $\xi= \xi_{k\eps^{-2}\delta}^{(\delta,+)}$  by adding $N_{k;+} $ particles all at the origin
and then removing the $N_{k;-} $ rightmost particles.

\end{definition}

\medskip
Thus  in the $\xi_t^{(\delta,\pm)}$ processes
  births and deaths are
concentrated at the times $k\eps^{-2}\delta$, in between
such times the particles are independent random walks. While the analysis of the true process $(\xi_t)_{t\ge 0}$
is rather complex due to the non local nature of
$L_a$, the study of the hydrodynamical limit for  $\xi_t^{(\delta,\pm)}$ is much simpler because
the number of rightmost particles to delete is macroscopic and becomes deterministic, the analysis
will be carried out in the next section.

\medskip
We shall often use in the sequel the following explicit realization of the random walk process $(n_t)_{t\ge 0}$.

\medskip

\begin{definition} [The probability space $(\Om_0,P_0)$]
\label{defin:anna3.1}  We set
 $\Om_0 = \{ \om_0= (\und t_0,\und \si_0)\}$, where
$\und t_0=(t_{1;0},\ t_{2;0}, \ldots)$, $\und \si_0=(\si_{1;0},\si_{2;0},...)$
are infinite sequences
of increasing positive ``times'' $t_{h;0}$ and
of   {symmetric ``jumps''},  $\si_{h;0} = \pm 1$.
$(\Om_0,P_0)$ is the  product
of a Poisson
 {process} of intensity  $2j\eps$  for
the increments of the time  {sequence} $\und t_0$  and of a
Bernoulli
 {process with parameter $1/2$} for
the {jump sequence} $\und \si_0$.

%

Given $n_0\in \mathbb N$ and $\om_0\in \Om_0$ we define $(n_t)_{ t\ge 0}$, iteratively: we set $n_t= n_{t_{h;0}}$
in the time interval $[t_{h;0},t_{h+1;0})$, $h\ge 0$, ($t_{0;0}\equiv 0$) and  define
\[
n_{t_{h+1;0}} = \begin{cases} n_{t_{h;0}} + \si_{h+1;0} &\text{if $n_{t_{h;0}} + \si_{h+1;0}\ge 0$}\\
0 &\text{if $n_{t_{h;0}} + \si_{h+1;0}< 0$}\end{cases}
\]
\end{definition}

\medskip

It is readily seen that the law of $(n_t)_{ t\ge 0}$ as a process on $(\Om_0,P_0)$
(for a given initial value $n_0$) is the same
as the Markov process of Theorem \ref{thme3.6a} and hence of the particles' number $|\xi_t|$
in our original
process once $n_0=|\xi_0|$.

\medskip
When realized on $(\Om_0,P_0)$, $N_{k, +}\equiv N_{k, +}(\om_0,n_0)$ ($n_0$ the initial particles' number)
is the number of times $t_{h;0}\in [k\eps^{-2}\delta,(k+1)\eps^{-2}\delta]$
where $\si_{h;0}=1$ (which does not depend on $n_0$) , while the number of times $t_{h;0}\in [k\eps^{-2}\delta,(k+1)\eps^{-2}\delta]$
where $\si_{h;0}=-1$ is an upper bound for $N_{k, -}\equiv N_{k, -}(\om_0,n_0)$ as the values $\si_{h;0}=-1$
do not produce a jump if $n_{t_{h;0}}=0$ (hence the dependence on $n_0$).

\medskip

Under the assumptions
on the initial datum $\xi$, see Definition \ref{defin0.1}, the process of
adding and removing particles becomes quite simple.  For any integer $k>0$ define  on $\Om_0$
	\begin{equation}
	\label{e3.1aa.1}
B^0_{k}(\om_0)  = \sum_{h} \mathbf 1_{\si_{h;0}=+1}
\,\,\mathbf 1_{t_{h;0}\in [k\eps^{-2}\delta,(k+1)\eps^{-2}\delta]}
	\end{equation}
	\begin{equation}
	\label{e3.2aa.2}
A^0_{k}(\om_0)  = \sum_{h}  \mathbf 1_{\si_{h;0}=-1}
\,\,\mathbf 1_{t_{h;0}\in [k\eps^{-2}\delta,(k+1)\eps^{-2}\delta]}
	\end{equation}
$B^0_{k}$ and $A^0_{k}$ are independent Poisson distributed variables
with average $\eps^{-1}j\delta$.

\medskip

\begin{definition}[{Good sets}]
Given $T>0$ and $\ga>0$ we define for any $\delta$ and $\eps$ positive
	\begin{equation}
	\label{e3.2aa.3}
\mathcal G  = \Big\{ \om_0 \in \Om_0: |A^0_k(\om_0) - \eps^{-1} j \delta| \le \eps^{-\frac 12-\ga};
|B^0_k(\om_0) - \eps^{-1} j \delta| \le \eps^{-\frac 12-\ga},\; k: k 	\delta \le T\Big\}
	\end{equation}

\end{definition}

\vskip.5cm

\begin{theorem}[{Reduction to Poisson variables}]
\label{thme3.6a.1}
Given $\xi$ as in  Definition \ref{defin0.1}, $T>0$ and $\ga>0$ there is $\delta^*>0$ so that for any
$\delta< \delta^*$ and any
$\eps>0$ small enough the following holds. For any $\om_0\in \mathcal G$ (see \eqref{e3.2aa.3}) and any $k$ such that
$k\delta \le T$,
	\begin{equation}
	\label{e3.2aa.4}
N_{k, -}(\om_0,|\xi|)= A^0_{k}(\om_0),\qquad N_{k,+}(\om_0,|\xi|)=B^0_{k}(\om_0)
	\end{equation}
where $N_{k, \pm}(\om_0,|\xi|))$  denote the variables
 $N_{k, \pm}$ when realized on $\Om_0$.

Finally, for any $n$ there is $c_n$ so that
	\begin{equation}
	\label{e3.2aa.5}
P_0[\mathcal G] \ge 1 - c_n\eps^n
	\end{equation}

\end{theorem}

\medskip

\noindent
{\bf Proof.} By  Definition \ref{defin0.1} the initial number of particles $|\xi|$ is bounded from below by
$\eps^{-1} \int \rho_{\rm init} - \eps^{-1+a}\ge \eps^{-1}C$, $C>0$.  We choose $\delta^*:= C/(2j)$ and shall prove by induction that for any $\delta<\delta^*$ and  all $\eps$ small enough we have in $\mathcal G$
\[
n_{t_k} \ge \eps^{-1}C - k2\eps^{-\frac 12 - \ga}, \;\;\; k\le \frac{ T}\delta,\qquad t_k=k\eps^{-2}\delta
\]
Suppose that the inequality holds for $k$ and let us prove it for $k+1$.  Since $N_{k, -}(\om_0,|\xi|)\le A^0_{k}(\om_0)$
\[
n_t \ge n_{t_k} - \eps^{-1} j \delta - \eps^{-\frac 12 -\ga}\ge \eps^{-1}(C- j\delta^*)
-(2k+1)\eps^{-\frac 12 -\ga}, \quad t \in [t_k,t_{k+1}]
\]
which is strictly positive for  any $k \le T/\delta$ if
$\eps$ is small enough.   Thus \eqref{e3.2aa.4} holds and
\[
n_{t_{k+1}} \ge {n_{t_k}} - A^0_{k}(\om_0) + B^0_{k}(\om_0)\ge
 n_{t_k} -2 \eps^{-\frac 12 -\ga}
\]
because $\om_0\in \mathcal G$.  This proves the induction hypothesis and for what
seen in the proof, \eqref{e3.2aa.4} holds as well.

The variables $A^0_{k}(\om_0)$ and $B^0_{k}(\om_0)$, $k\le T/\delta$, are
independent Poisson variables with mean $\eps^{-1}j\delta$ hence
\eqref{e3.2aa.5}.  \qed

\vskip.5cm

Once restricted to $\mathcal G$ the processes $\xi_t^{(\delta,\pm)}$,
$0\le t \le \eps^{-2}T$, become quite simple.
The particles move as independent random walks  in
the finitely many intervals $[k\eps^{-2}\delta,(k+1)\eps^{-2}\delta]$, while births
and deaths  at the times $k\eps^{-2}\delta$ are ``essentially deterministic'' like in the
$\delta$-approximated evolutions $S_t^{(\delta,\pm)}$
of  Definition \ref{defin:e2.3}.
Such considerations are made precise
in Section \ref{quarta} where
we prove convergence of  $\xi_t^{(\delta,\pm)}$  to
$S_t^{(\delta,\pm)}(\rho_{\rm init})$ in the hydrodynamic limit.

\vskip2cm

\section{Hydrodynamic limit  for the approximating processes }
\label{quarta}

The main result in this section is in  Theorem \ref{teo:3.1} below.
It states that the
$\delta$-approximate processes $\xi_t^{(\delta,\pm)}$
of Definition \ref{define3.1aa} converge in
the hydrodynamic limit  to
the evolutions
$S_t^{(\delta,\pm)}(\cdot)$ of
Definition \ref{defin:e2.3}.

Here we exploit duality to prove
convergence in a very strong form of  the
independent system to the heat equation.

For any fixed $\delta$ and $T>0$,
the processes  $\xi_t^{(\delta,\pm)}$, $t\le \eps^{-2}T$ are
obtained by alternating independent random walk evolutions to cut and paste
operations.  The latter involve macroscopic quantities and can be controlled by means of Theorem \ref{thme3.6a.1}
once we have the hydrodynamic limit for
the independent process.  This is well studied and very detailed estimates are
available but in the present case we have the extra difficulty  that
the initial configurations may have a macroscopic occupation number at the origin
$\xi(0)\approx \eps^{-1}$.  This is because in the the cut and paste we actually paste
$\approx j\delta \eps^{-1}$ particles at the origin.
This is not a case studied in the literature
(as far as we know) and indeed it affects greatly the decay of correlations
in the hydrodynamic limit.

As  in our iterative procedure we have initial data with macroscopic occupation
at the origin, we may as well take more general initial conditions
(than those in Definition \ref{defin0.1})  with macroscopic occupation
at the origin, this will be actually useful in the sequel.  Thus the ``macroscopic initial profile
$v_0$ '' is here taken in  $\mathcal U$,
namely it is the sum of  a non negative $L^\infty$ function plus $c D_0$, with
$c$ either equal to $0$ or to $ j\delta$, we suppose that $\int v_0 = F(0;v_0)>0$.
Analogously to \eqref{3.3.1} for any $\eps>0$ we choose the initial configuration
$\xi_0$ so that
 \begin{equation}
  \label{3.3.1.4k}
\max_{x\in [0,\eps^{-1}-\ell+1]}\Big| \mathcal A_\ell(x,\xi_0)-
  \mathcal A'_\ell(x,v_0)\Big| \le \eps^a
\end{equation}

\medskip

\begin{theorem}
\label{teo:3.1}

Given any  $T>0$ for any $\delta>0$ small enough,
any $k: k\delta\le T$
and any $\zeta>0$
   \begin{equation}
   \label{5.15}
\lim_{\eps \to 0}P^{(\eps)}_{\xi_0}\Big[\max_{x\in [0,\eps^{-1}]}|\eps F_\eps(x;\xi^{(\delta,\pm)}_{k\eps^{-2}\delta}) - F(\eps x;S_{k\delta}^{(\delta,\pm)}(v_{0}))| \le \zeta  \Big] = 1
    \end{equation}
where $v_0$ and $\xi_0$ are as above;
$P^{(\eps)}_{\xi_0}$
as in Definition \ref{defin0.1}; $F$ and $F_\eps$ as in \eqref{ee8.1}.

\end{theorem}

 \medskip

The  theorem is proved at
the end of the section, as we shall see  stronger results actually hold
but what  stated is what needed for
Theorem \ref{Teo:Hydro}.
In the course of the proof we shall introduce several
positive
parameters: $b$, $a$, $a^*$, $\ga$: $b$ should be close to 1 and the others close to $0$,
for the sake of definiteness we take:
 \begin{equation}
  \label{const}
a=\ga = \frac 1{20},\;\; b=\frac 9{10},\;\;a^*=\frac 1{100}
	\end{equation}
We prove the theorem only for
the process $\xi_t^{(\delta,-)}$, the analysis
of $\xi_t^{(\delta,+)}$ is similar and omitted.
The first step is a spatial discretization of the flow $S_{k\delta}^{(\delta,-)}$:

\medskip

\begin{definition}[{The discrete evolution}]
\label{definanna4.1}
Denote by  $p^0_{t}(x,y)$, $t\ge 0$, $x,y\in [0,\eps^{-1}]$,
the transition probability of a
continuous time, simple symmetric random walk with reflections at
$0$ and $\eps^{-1}$
(i.e.\ the random walker jumps by $\pm 1$ with equal probability
after an exponential time of mean 1,
the jumps which would lead outside $[0,\eps^{-1}]$ are
suppressed).  For $\delta$ small enough
we define functions $u_k(x)$, $x\in [0,\eps^{-1}]\cap \mathbb Z$, with the property that
mass is conserved: $F_\eps(0;u_{k})=F_\eps(0;u_0)$ for all $k$.
The definition is iterative, we set $u_0(x) := v_0(\eps x)$;
then supposing that $u_{k-1}$ has been already defined
and  that $F_\eps(0;u_{k-1})=F_\eps(0;u_0)$ we define $u_k$ as follows. We first call
 \begin{equation}
  \label{anna4.1}
u^0_k(x) = \sum_y p(x,y) u_{k-1}(y), \quad \quad  p(x,y):= p^0_{\eps^{-2}\delta}(x,y)
	\end{equation}
$u_k$ is then
obtained from $u^0_k$ by adding particles at 0 and removing particles on the right.
To make this precise let
$R_k$  be an integer
such that
$F_\eps(R_k;u^0_k) \ge \eps^{-1}j \delta$ while $F_\eps(R_k+1;u^0_k) < \eps^{-1}j \delta$.
The existence of
$R_k$
for $\delta$ small enough follows from the assumption
$F_\eps(0;u_0) \ge c\eps^{-1}$, $c>0$,
observing that $F_\eps(0;u^0_k)=F_\eps(0;u_{k-1})=F_\eps(0;u_0)$
by the inductive assumption
and $F_\eps(0;u_0)= \eps^{-1}F(0;v_0)$.
We then set
$v_k(x)= u^0_k(x)$ for $x<R_k$, $v_k(x)=0$ for $x>R_k$ and
\begin{equation*}
  v_k(R_k) := F_\eps (R_k;u^0_k) - \eps^{-1}j \delta
	\end{equation*}
We can then finally define $u_k$ as
\begin{equation}
  \label{anna4.1.0}
 u_k := v_k+ \eps^{-1} j\delta \mathbf 1_0
	\end{equation}
where $\mathbf 1_0$ is the Kr\"onecker delta at 0.
To complete the induction we observe that
$F_\eps(0;u_k)= \eps^{-1} j\delta + F_\eps(0;v_k)$, $F_\eps(0;v_k)= F_\eps(0;u^0_k)- \eps^{-1}j \delta$
so that $F_\eps(0;u_k)=F_\eps(0;u^0_k)=F_\eps(0;u_{k-1})$.

\end{definition}

\medskip

In the next proposition we show that in \eqref{5.15}  we can replace $S_{k\delta}^{(\delta,-)}(v_{0})$
by the sequence $u_k$ with a negligible error:

\medskip

\begin{proposition}
In the same context as in Theorem \ref{teo:3.1},
   \begin{equation}
   \label{5.15.1}
\lim_{\eps \to 0}\max_{x\in [0,\eps^{-1}]}|\eps F_\eps(x;u_k) - F(\eps x;S_{k\delta}^{(\delta,-)}(v_{0}))| =0
    \end{equation}
\end{proposition}

\medskip
\noindent
{\bf Proof.}
In this proof we shorthand by $g(r,r')$   the Green function
$G_\delta^{\rm neum}(r,r')$, $r,r' \in [0,1]$, defined in \eqref{e4.3} and
 also write for brevity $p(x,y):= p^0_{\eps^{-2}\delta}(x,y)$,
 as in Definition \ref{definanna4.1}.
Let $u_k$, $u^0_k$ and $R_k$ be as in Definition \ref{definanna4.1}.  We define for any real $r$ between $0$ and $\eps^{-1}$,
 \begin{equation*}
\psi_k(r) := [S_{k\delta}^{(\delta,-)}(v_{0}) - j\delta D_0](\eps r)
	\end{equation*}
Analogously to \eqref{3.3.1.00} we denote by $R'_k$ the real number in $[0,\eps^{-1}]$
such that $\psi_k(r)>0$ for $r<R'_k$ and $\psi_k(r)=0$ for $r>R'_k$.  We also call
 \[
 \psi^0_k(r)= j\delta g(\eps r,0)+ \int_0^1 dr'
 g(\eps r,r') \psi_{k-1}(\eps^{-1} r') 
 \]
so that
  \[
 \psi^0_k(r) = \psi_k(r),\; r<R'_k ;\quad \int_{R'_k}^{\eps^{-1}} \psi^0_k( r) = \eps^{-1} j\delta
  \]

  \medskip

\noindent
{\bf Claim.} There are strictly positive  constants $C_{\pm}$ which depend on $\delta$ so that for all $k$,
 \begin{eqnarray}
\label{e4anna.1}
&& C_- \le \psi^0_k \le C_+,\quad C_- \le u^0_k \le C_+; \quad |\frac{d}{dr}\psi^0_k| \le \eps C_+\nn
\\&&
|\sum_{x\in \mathbb Z: x\in [ R_k,\eps^{-1}]} \psi^0_{k}(x) - \eps^{-1}j\delta| \le  C_+,\quad
|F_\eps(x; \psi^0_{k}) - \int_{x}^{\eps^{-1}}\psi^0_{k}| \le  C_+
   \end{eqnarray}

  \medskip
The proof of the claim
follows from classical estimates on random walks and Green functions:
 \begin{eqnarray}
 \label{a-e4.8}
&& \frac{c_1}{\sqrt \delta} \le g(r,r') \le \frac{c_2}{\sqrt \delta};\quad
\frac{c_1\eps}{\sqrt \delta} \le p(x,y) \le \frac{c_2\eps }{\sqrt \delta}
\\&& |\frac{d}{dr}g(r,r')| \le \frac{c_3}{\delta};\quad
|p(x,y)-p(x,y+1)| \le (\frac{c_3 \eps}{\sqrt\delta})^2\nn
\\&&
j\delta+\eps\int_0^{\eps^{-1}} \psi_k = F(0;v_0) \nn
   \end{eqnarray}

  \medskip

\noindent
The crucial step in the proof of the proposition is the following statement:

\medskip

 \begin{equation}
\label{e4anna.2}
\text{There
are $\alpha >\beta>1$ so that}\;\; |u^0_k(x) -\psi^0_k(x)| \le \frac{\eps}{\sqrt \delta} \alpha^k, 
\quad
 |R_k-R'_k| \le \beta \alpha^k
 	\end{equation}

 \medskip

\noindent
We  prove \eqref{e4anna.2} by
induction.  We thus suppose that it holds for $k-1$.
Calling $R^*_k$ the largest integer smaller or equal than $R_k$ and $R'_k$
   \begin{eqnarray*}
&&|u^0_k(x)-\psi^0_k(x)|  \le j\delta  |  \eps^{-1}p(x,0) - g(\eps x,0)|
\\&& \hskip1cm+ \sum_{y \le R^*_{k-1}}\Big|p(x,y) u^0_{k-1}(y) - \int_{\eps y}^ {\eps(y+1)} g(\eps x,r) \psi^0_{k-1}(\eps^{-1}r)\Big|
\\&& \hskip1cm +\sum_{y = R^*_{k-1}+1}^{R_{k-1}} p(x,y) u^0_{k-1}(y)
+\int_{\eps R^*_{k-1}}^{\eps R'_{k-1}} g(\eps x,r) \psi^0_{k-1}(\eps^{-1}r)
   \end{eqnarray*}
We use the local central limit theorem to bound:
 \begin{equation}
\label{e4anna.2.0}
\Big|p(x,y)- \eps g(\eps x,\eps y)\Big| \le \frac{c_5\eps^2}{\delta}
 	\end{equation}
 Thus
    \begin{eqnarray*}
&&|u^0_k(x)-\psi^0_k(x)|  \le j\delta  \frac{c_5\eps}{\delta} + \max_{x} |u^0_{k-1}(x)-\psi^0_{k-1}(x)|
\\&& \hskip1cm+ \sum_{y \le R^*_{k-1}}\Big|p(x,y) \psi^0_{k-1}(y) - \int_{\eps y}^ {\eps(y+1)} g(\eps x,r) \psi^0_{k-1}(\eps^{-1}r)\Big|
\\&& \hskip1cm + 2\frac{c_2\eps }{\sqrt \delta} |R_{k-1}-R'_{k-1}| C_+
   \end{eqnarray*}
   We write
       \begin{eqnarray*}
&&\Big|\int_{\eps y}^ {\eps(y+1)} g(\eps x,r) \psi^0_{k-1}(\eps^{-1}r)
-\eps g(\eps x,\eps y)\psi^0_{k-1}(y) \Big|
\le  c_6 \eps^2
   \end{eqnarray*}
   and get using the induction assumption
       \begin{eqnarray*}
&&|u^0_k(x)-\psi^0_k(x)|  \le j c_5\eps  + \frac{\eps}{\sqrt \delta} \alpha^{k-1} + 2\frac{c_2\eps }{\sqrt \delta} \beta \alpha^{k-1} C_+
\\&& \hskip2cm+ \eps^{-1}\{c_6 \eps^2 + C_+  \frac{c_5\eps^2}{\delta}\}
   \end{eqnarray*}
Choosing  $\alpha\ge 1+ j\sqrt \delta  c_5
  + 2 c_2 C_+ \beta
 +  \sqrt \delta \{c_6  + C_+  \frac{c_5 }{\delta}\}$, we have
    \begin{eqnarray*}
&&|u^0_k(x)-\psi^0_k(x)|  \le \frac{\eps}{\sqrt \delta}\alpha^{k-1}\Big( j\sqrt \delta  c_5  + 1
  + 2 c_2 C_+ \beta
+  \sqrt \delta \{c_6  + C_+  \frac{c_5 }{\delta}\}\Big)  \le \frac{\eps}{\sqrt \delta}\alpha^{k}
   \end{eqnarray*}
As a consequence:
 \begin{equation}
  \label{e44anna4.1}
|F_\eps( x;u^0_k) -F_\eps( x;\psi^0_k)| \le (\eps^{-1}-x+1) \frac{\eps}{\sqrt \delta} \alpha^k
 	\end{equation}
Recalling that $\psi^0_k(x)\ge C_-$ and $u^0_k(x) \ge C_-$  we get
 \begin{eqnarray*}
&&|F_\eps(R'_k;u^0_k) -j\delta \eps^{-1}| \le C_+ +(\eps^{-1}-R'_k+1) \frac{\eps}{\sqrt \delta} \alpha^k\\
&&|F_\eps(R'_k;u^0_k)-F_\eps(R_k;u^0_k)| \le 2C_+ +  \eps^{-1}\frac{\eps}{\sqrt \delta} \alpha^k\\&&
C_- |R'_k-R_k|\le |F_\eps(R'_k;u^0_k)-F_\eps(R_k;u^0_k)| \le 2C_+ + \frac{  \alpha^k}{\sqrt \delta}
 	\end{eqnarray*}
which is smaller than $\beta \alpha^k$ if $\beta\ge C_-^{-1}(2C_+ +  \delta^{-1/2})$, thus completing the proof
of \eqref{e4anna.2}.

Using \eqref{e44anna4.1}  we then conclude the proof of the proposition, details are omitted.

\qed

\vskip.5cm

The proof of
Theorem \ref{teo:3.1}  is thus reduced to showing that: for all $n$ so that $n\delta\le T$,
   \begin{equation}
   \label{5.15.2}
\lim_{\eps \to 0}P^{( \eps)}_\xi\Big[\eps |F_\eps(x;\xi^{(\delta,-)}_{n\eps^{-2}\delta}) - F_\eps( x;u_n)| \le \zeta \;\;\text{for all $x\in [0,\eps^{-1}]$} \Big] = 1 
    \end{equation}
which   will be done in the sequel.  Both sequences $\{\xi^{(\delta,-)}_{n\eps^{-2}\delta}\}$
and   $\{u_n\}$ are determined by alternating free evolution and a cut and paste procedure.
We first study the free evolution part proving that the independent random walk configuration
$\xi^0_{\eps^{-2}\delta}$ is well approximated by its average.
Call $P_\xi$ and $E_\xi$ law and expectation
of the independent process
starting from  $\xi$, define for $x\in [0,\eps^{-1}]$
 	\begin{equation}
	\label{5.11}
w(x|\xi):=\mathbb E_\xi[\xi^0_{\eps^{-2}\delta}(x)]=
\sum_{y=0}^{\eps^{-1}} p(x,y) \xi(y),\quad  p(x,y):= p^0_{\eps^{-2}\delta}(x,y)
	\end{equation}
with $p^0_t$ the transition probability used in Definition \ref{definanna4.1}.

\medskip

		\begin{proposition}
			\label{teo:3.2}
Let $c^*$ and $a^*$ be strictly  positive and
	\begin{equation}
	\label{a-e4.13}
\mathcal X_{c^*,a^*}:= \Big\{\xi: |\xi|\le c^*\eps^{-1},\;
\max_{x\ne 0}\xi(x)\le \eps^{-a^*}\Big\}
	\end{equation}
Then for any  $\xi\in \mathcal X_{c^*,a^*}$
	\begin{equation}
	\label{a-e4.13.0}
\max_{x \in [0,\eps^{-1}]}\; w(x|\xi) \le \frac{c_2 c^*}{\sqrt \delta}
	\end{equation}
($c_2$ as in \eqref{a-e4.8}).  Moreover
let
$ c^*$, $a^*$ and $b$ be strictly positive and
such that
	\begin{equation}
	\label{a*}
   a^*< \frac b2,\qquad b+a^*<1
	\end{equation}
(a condition which is satisfied by the choice \eqref{const}).
Let $\ell$ be the integer part of $\eps^{-b}$ and $\mathcal A_\ell$ be as in \eqref{3.3.0}, then
for any integer $n$ there is $c'_n$ so that
	\begin{equation}
	\label{a-e4.13}
P_\xi\Big[\xi^0_{\eps^{-2}\delta}\in  \mathcal X_{c^*,a^*}\Big] \ge 1-  c'_n \eps^n
	\end{equation}
Finally there is a constant $c$ so that
	\begin{equation}
	\label{a4.7}
	\sup_{x \le \eps^{-1}-\ell+1} E_{\xi}\big[ \big|\mathcal A_\ell(x, \xi^0_{\eps^{-2}\delta})-
  \mathcal A_\ell(x, w(\cdot|\xi))\big|^{4}\big]\le c \eps^{2b}
  	\end{equation}

\end{proposition}

	\medskip

\noindent {\bf Proof.}  For brevity in this proof we shall write $w(x)$ instead of
$w(x|\xi)$. Recalling that $p(x,y)$ is defined in \eqref{anna4.1} and bounded in
\eqref{a-e4.8}, we have for any  $\xi\in \mathcal X_{c^*,a^*}$
	\begin{eqnarray}
w(x)&=&\sum_yp(x,y)\xi(y)\le  \frac{c_2 \eps}{\sqrt \delta}\sum_y \xi(y)
\le  \frac{c_2 \eps}{\sqrt \delta} c^*\eps^{-1}
	\label{a4.6}
	\end{eqnarray}
hence \eqref{a-e4.13.0}.
The proof of \eqref{a-e4.13} and \eqref{a4.7} uses in a crucial way duality:

	\medskip
	
{\em Duality}.
Given $\xi\in \mathbb N^{[0,N]}$ and a labeled configuration $\und x=(x_1,..,x_n)$, $n\ge 1$, $x_i\in [0,\eps^{-1}]$,
we define
	\begin{eqnarray}
	\label{a4.9}
&&\mathcal D(\xi,\und x)=\prod_x
d_{\und x(x)}(\xi(x)),\quad d_k(m)=m(m-1)\cdots (m-k+1),\; d_0(m)=1\\&&  \und x(x)
=\sum_{i=1}^n \mathbf 1_{x_i=x} \nn
	\end{eqnarray}
$d_k(m)$ are called
Poisson polynomials.  We then have:
	\begin{equation}
	\label{a4.10}
E_{\xi}\big[\mathcal D(\xi^0_t,\und x)\big]=E_{\und x}\big[\mathcal D(\xi,\und x^0_t)\big]
	\end{equation}
where $\und x^0_t$ is the independent random walks evolution.

	\medskip

\noindent
$\bullet$ {\it Proof of \eqref{a-e4.13}.}  Call
$\und x= (x_1,..,x_{2k})$ with $x_i=x$ for all $i=1,..,2k$.  Then by \eqref{a4.10}
and \eqref{a4.6}
	 \begin{eqnarray}
 	\label{a4.14aa}
E_{\xi} [d_{2k}(\xi^0_{\eps^{-2}\delta}(x))] &=& E_{\und x} [\prod _x d_{\und x^0_{\eps^{-2}\delta}(x)}(\xi_0(x))]\le
E_{\und x} [\prod _x \xi(x)^{\und x^0_{\eps^{-2}\delta}(x)}]\nn\\& = &\Big[\sum_yp(x,y)\xi(y)\Big]^{2k} \le
\Big(\frac{c_2 \eps |\xi| }{\sqrt \delta} \Big)^{2k}\le \Big(\frac{c_2 c^* }{\sqrt \delta} \Big)^{2k}
	\end{eqnarray}
 By \eqref{a4.14aa} we have that
for any $k$ there is $c''_k$ (independent of $\eps$) so that
	  \begin{equation}
	\label{bound}
\max_{x\in[0,\eps^{-1}]}
E_{\xi}\big[\xi^0_{\eps^{-2}\delta}(x)^k \big]\le  c''_k
	\end{equation}
Moreover by the Chebishev inequality and \eqref{a4.14aa}
	  \begin{equation}
	\label{a4.14b}
P_\xi\big[\max_{x\in[0,\eps^{-1}]}\xi^0_{\eps^{-2}\delta}(x)\le \eps^{-a^*}\big]\ge 1-c'_m\eps^m
	\end{equation}
which proves  \eqref{a-e4.13} because $|\xi^0_{\eps^{-2}\delta}|=|\xi| \le \eps^{-1}c^*$.

\medskip

To prove \eqref{a4.7} we shall use again duality but also several
maybe non totally straightforward
algebraic manipulations.
We start by expanding the product in the expectation:
	\begin{eqnarray}
  	\label{a4.8}
&&\hskip-2cm E_{\xi}\Big[ \big|\mathcal A_\ell(x, \xi^0_{\eps^{-2}\delta})-\mathcal A_\ell(x,w)\big|^{4}\Big]
=\frac 1{\ell^{4}}\sum_{\und x \in \mathcal B_\ell}E_{\xi}\Big[\prod_{i=1}^4\big(\xi^0_{\eps^{-2}\delta}(x_i)-w(x_i)\big)\Big]
	\end{eqnarray}
where  $x \in [0,\eps^{-1}-\ell+1]$ and
$\mathcal B_\ell=\{ \und x =(x_1,...x_{4}): x_i\in [x,x+\ell-1], i=1,..,4\}$.

Call $\mathcal B_\ell^{(i)}$, $i=1,2,3,4$, the set of $\und x \in \mathcal B_\ell$ such that there are $i$ mutually distinct sites.  We then have for $i\le 2$:
	\begin{eqnarray}
  	\label{a4.8bis}
&&\hskip-2cm \frac 1{\ell^{4}}\sum_{\und x\in \mathcal B^{(i)}_\ell}|E_{\xi}\Big[\prod_{i=1}^4\big(\xi^0_{\eps^{-2}\delta}(x_i)-w(x_i)\big)\Big]| \le c \ell^{-2}
	\end{eqnarray}
as the expectation of products of $\xi^0_{\eps^{-2}\delta}(\cdot)$ is bounded, which is proved using
\eqref{bound}.

We are thus left with the sum over
$\und x\in \mathcal B^{(i)}_\ell$ with $i=3,4$.  When   $i=4$, $\und x=(x_1,..,x_4)$ with the entries mutually distinct.
Call $\si=(\si_1,..,\si_4)$, $\si_i\in \{-1,1\}$, and $|\si|_-$ the number of $-1$ in $\si$, then
	\begin{eqnarray}
  	\label{a4.8.n1}
&&\prod_{i=1}^4\big(\xi^0_{\eps^{-2}\delta}(x_i)-w(x_i)\big)= \sum_\si (-1)^{|\si|_-} \mathcal D(\xi^0_{\eps^{-2}\delta};\{x_i: \si_i=1\}) \prod_{j:\si_j=-1} w(x_j)
	\end{eqnarray}
and using duality:
	\begin{eqnarray}
  	\label{a4.8.n2}
&&
E_\xi[\prod_{i=1}^4\big(\xi^0_{\eps^{-2}\delta}(x_i)-w(x_i)\big)]=\sum_{\und y} p(\und x,\und y)\sum_\si (-1)^{|\si|_-} \mathcal D(\xi;\{y_i:\si_i=1\})\nn
\\&& \hskip5cm \times \Pi(\xi;\{y_j:\si_j=-1\})\nn\\&&
 \Pi(\xi;\{y_j:\si_j=-1\}):= \prod_{j:\si_j=-1} \xi(y_j)
	\end{eqnarray}
Suppose there is a singleton $h$, namely such that $y_h\ne y_j$ for all $j\ne h$, then
	\begin{eqnarray}
  	\label{a4.8.n3}
&&\hskip-2cm
 \sum_\si (-1)^{|\si|_-} \mathcal D(\xi^0_{\eps^{-2}\delta};\{y_i:\si_i=1\}) \Pi(\xi;\{y_j: \si_j=-1\}) =0
	\end{eqnarray}
Indeed let $\si$ a sequence with $\si_h=1$ and $\si'$ the one obtained from $\si$ by changing only $\si_h$, then
	\begin{eqnarray*}
&&
 (-1)^{|\si|_-} D(\xi^0_{\eps^{-2}\delta};\{y_i: \si_i=1\})\Pi(\xi;\{y_j: \si_j=-1\}) \\&& \hskip1cm =
(-1)^{|\si|_-}  D(\xi^0_{\eps^{-2}\delta};\{y_i: \si_i=1, i\ne h\}) \Pi(\xi;\{y_h, y_j:\si_j=-1\})
%
%
\\&& \hskip1cm =-(-1)^{|\si'|_-}  \mathcal D(\xi^0_{\eps^{-2}\delta};\{y_i: \si'_i=1\}) \Pi(\xi;\{y_j:\si'_j=-1\})
 	\end{eqnarray*}
We have thus proved that calling $\mathcal
X_{\rm n.s.}$ the set of all $\und y$ with no singletons then
 	\begin{eqnarray}
  	\label{a4.8.n4}
&&
E_\xi[\prod_{i=1}^4\big(\xi^0_{\eps^{-2}\delta}(x_i)-w(x_i)\big)]= \Phi_4(\und x)\\&&
\Phi_4(\und x)=\sum_{\und y\in \mathcal X_{\rm n.s.}} p(\und x,\und y)
\sum_\si (-1)^{|\si|_-} \mathcal D(\xi;\{y_i:\si_i=1\}) \Pi(\xi;\{y_j:\si_j=-1\})\nn
%
	\end{eqnarray}
 A similar property holds also when $\und x \in \mathcal
 B^{(3,*)}_\ell$ which is the set of all
$\und x$ such that $x_1=x_2$, $x_3\ne x_4,x_1$ and $x_4\ne x_1$ (modulo permutation of labels
all $\und x\in \mathcal
 B^{(3)}_\ell$ are in $\mathcal
 B^{(3,*)}_\ell$).  We write
\[
\Big( \xi(x)-w(x)\Big)^2= \{\xi(x) [\xi(x)-1] - 2w(x)\xi(x) +w(x)^2\} + \{\xi(x) -w(x)\}+w(x)
\]
Then analogously to \eqref{a4.8.n1}	but with $\und x\in \mathcal B^{(3,*)}_\ell$,
    \begin{eqnarray}
  	\nn
&&\hskip-1.5cm\prod_{i=1}^4\big(\xi^0_{\eps^{-2}\delta}(x_i)-w(x_i)\big)= \sum_{\si \in \{-1,1\}^4} (-1)^{|\si|_-} D(\xi^0_{\eps^{-2}\delta};\{x_i, \si_i=1\}) \prod_{j:\si_j=-1} w(x_j) \nn
\\&& \hskip1cm +\sum_{\si =(\si_2,\si_3,\si_4)} (-1)^{|\si|_-} D(\xi^0_{\eps^{-2}\delta};\{x_i, \si_i=1, i\ge 2\})
\prod_{j\ge 2:\si_j=-1} w(x_j)
\nn \\&& \hskip1cm + w(x_1) \sum_{\si =(\si_3,\si_4)} (-1)^{|\si|_-} D(\xi^0_{\eps^{-2}\delta};\{x_i, \si_i=1, i\ge 3\}) \prod_{j\ge 3:\si_j=-1} w(x_j)\nn\\
\label{a4.8.n5}
	\end{eqnarray}
	\begin{eqnarray}
  	\label{a4.8.n6}
E_\xi[\prod_{i=1}^4\big(\xi^0_{\eps^{-2}\delta}(x_i)-w(x_i)\big)]= \Phi_4(\und x) + \Phi_3(x_2,x_3,x_4) +w(x_1) \Phi_2(x_3,x_4)\end{eqnarray}
where
\begin{eqnarray*}
&&
\hskip-.6cm \Phi_3(x_2,x_3,x_4)= \sum_{(y_2,y_3,y_4)\in \mathcal X_{\rm n.s.}} \prod_{i=2}^4 p( x_i, y_i)
\sum_{\si_2,\si_3,\si_4} (-1)^{|\si|_-}\mathcal D(\xi;\{y_i, \si_i=1, i\ge 2\})\nn\\&&\hskip5cm \times
\Pi(\xi;j\ge 2:\si_j=-1) \nn
\\&&
\hskip-.6cm\Phi_2(x_3,x_4)=  \sum_{(y_3,y_4)\in \mathcal X_{\rm n.s.}}
\prod_{i=3}^4 p( x_i, y_i)\sum_{\si =(\si_3,\si_4)} (-1)^{|\si|_-} D(\xi;\{y_i, \si_i=1, i\ge 3\})\nn\\&&\hskip5cm \times \Pi(\xi;j\ge 3:\si_j=-1)\nn
	\end{eqnarray*}
with $\Phi_4(\und x)$ as in \eqref{a4.8.n4}.

Going back to \eqref{a4.8},  using
\eqref{a4.8bis} and \eqref{a-e4.13.0}
	\begin{eqnarray}
  	\label{a4.8n7}
&&\hskip-2cm E_{\xi}\Big[ \big|\mathcal A_\ell(x, \xi^0_{\eps^{-2}\delta})-\mathcal A_\ell(x,w)\big|^{4}\Big]
\le  \frac{c }{\ell^{2}} + \max_{\und x \in \mathcal B_\ell^{(4)}} \Phi_4(\und x)
+\frac{6}{\ell}\Big(
|\max_{\und x \in \mathcal B_\ell^{(3,*)}}\Phi_4(\und x)|\nn \\&&+| \max_{(x_2,x_3,x_4): {\rm distinct}}
\Phi_3(x_2,x_3,x_4)| + \frac{c_2 c^*}{\sqrt \delta} |\max_{(x_3,x_4): {\rm distinct}} \Phi_2(x_3,x_4)|\Big)
	\end{eqnarray}
Let us bound one by one the functions $\Phi_i$ starting from $\Phi_4$.  Recalling \eqref{a4.8.n4}
the condition $\und y\in \mathcal X_{\rm n.s.}$ is realized (modulo label permutations) in only
two cases: (i)\; $y_1=y_2\ne y_3=y_4$; (ii) $y_1=...=y_4$.
 	\begin{eqnarray}
  	\label{a4.8.n8}
\sum_\si (-1)^{|\si|_-} D(\xi;\{y_i:\si_i=1\}) \Pi(\xi;\{y_j:\si_j=-1\})=
\begin{cases} \xi(y_1)\xi(y_3),&\text{in case (i)}\\ 3\xi(y_1)^2 - 6\xi(y_1)
,&\text{in case (ii)} \end{cases}
	\end{eqnarray}
so that from \eqref{a-e4.8} and since $\xi\in  \mathcal X_{c^*,a^*}$
 	\begin{eqnarray}
  	\label{a4.8.n9}
&&
|\Phi_4(\und x)| \le  (\frac{c_2\eps}{\sqrt \delta})^4 \Big(  6 ( c^*\eps^{-1} )^2 +
3( c^*\eps^{-1} )^2 \Big)\le c \eps^2
%
	\end{eqnarray}
The condition $(y_2,y_3,y_4)\in \mathcal X_{\rm n.s.}$ in $\Phi_3$ implies $y_2=y_3=y_4$
and for such a $\und y$:
 	\begin{eqnarray}
  	\label{a4.8.n10}
\sum_{\si =(\si_2,\si_3,\si_4)} (-1)^{|\si|_-} D(\xi;\{y_i, \si_i=1, i\ge 2\})
\Pi(\xi;j\ge 2:\si_j=-1)= 2\xi(y_2)
	\end{eqnarray}
so that from \eqref{a-e4.8} and since $\xi\in  \mathcal X_{c^*,a^*}$
 	\begin{eqnarray}
  	\label{a4.8.n11}
&&
|\Phi_3(x_2,x_3,x_4)| \le  (\frac{c_2\eps}{\sqrt \delta})^3 2( c^*\eps^{-1} ) \le c \eps^2
%
	\end{eqnarray}
Finally if $(y_3,y_4)\in \mathcal X_{\rm n.s.}$ then $y_3=y_4$ and for such a $\und y$,
 	\begin{eqnarray}
  	\label{a4.8.n12}
\sum_{\si =(\si_3,\si_4)} (-1)^{|\si|_-} D(\xi;\{y_i, \si_i=1, i\ge 3\}) \Pi(\xi;j\ge 3:\si_j=-1)
= -\xi(y_3) \le 0
	\end{eqnarray}
Thus \eqref{a4.7} follows from \eqref{a4.8n7} together with the above inequalities.

\qed

\vskip.5cm

The cut and paste sequence of operations which appear in the definition of $\{\xi_{t_k}^{(\delta,-)}$,
$k \le k^*\}$, $k^*$ the largest integer such that $\delta k^* \le T$, $t_k=k\eps^{-2}\delta$, is independent
of the motion of the particles so that we have a rather explicit expression for the law of the variables
$\{\xi_{t_k}^{(\delta,-)}$,
$k \le k^*\}$, see \eqref{anna4.5} below.
We first write (with $\xi_0$ below
the initial condition in Theorem \ref{teo:3.1})
\begin{equation}
  \label{anna4.2}
p(\{n_k^{\pm}, k=1,.., k^*\}) = P_{\xi_0}^{(\eps)}\Big[ N_{k-1,-} =
n^-_k,  N_{k-1,+} = n^+_k, k \le k^*\Big]
	\end{equation}
where 
$N_{k,\pm}$  are  defined  in \eqref{e3.2aa} and \eqref{e3.1aa},
their law depends only on $|\xi_0|$.

We also write
\begin{equation}
  \label{anna4.3}
\pi(\xi'|\xi) = P_\xi\Big[\xi^0_{\eps^{-2}\delta} = \xi'\Big],\quad  |\xi|=|\xi'|, \;\; {\rm a.s.}
	\end{equation}
($\xi^0_t$  the independent random walk process).  We finally denote by
$K^{(n^-,n^+)} \xi$ the configuration obtained from $\xi$ by adding $n^+$
particles at 0 and then removing the $n^-$ rightmost particles (the definition
requires that $|\xi| + n^+-n^-\ge 0$, condition  automatically satisfied
below as the variables $n^{\pm}$ are the increments of the particles' number $n_t$).  Then, writing
\begin{eqnarray}
 \nn
&& P\Big[\{n_k^{\pm}, \xi_k^0,\; k=1,.., k^*\}\Big] = p(\{n_k^{\pm}, k\le k^*\}) \prod_{k=1}^{k^*}
\pi(\xi_k^0|\xi_{k-1})\\&& \xi_{k}:=K^{(n^-_{k},n^+_{k})}\xi^0_{k} \label{anna4.4}
	\end{eqnarray}
with  $n_0^{\pm}:=0$, we have
   \begin{equation}
   \label{anna4.5}
P^{(\eps)}_{\xi_0}\Big[\{\xi^{(\delta,-)}_{k\eps^{-2}\delta}) = \bar\xi_k,\; k=1,.., k^*\} \Big] =  \sum_{n_k^{\pm}, \xi_k^0,\; k=1,.., k^*} \mathbf 1 _{ \xi_{k}=\bar \xi_k,\; k=1,..,k^*}  P\Big[\{n_k^{\pm}, \xi_k^0,\; k\le k^*\}\Big]
    \end{equation}
By \eqref{e3.2aa.5} for any $n$ there is $c_n$ so that
\begin{eqnarray}
  \label{anna4.6}
&&
\sum_{\{n_k^{\pm},\; k=1,.., k^*\} \in \mathcal G}  p(\{n_k^{\pm}, k\le k^*\}) \ge 1 - c_n \eps^n\\&&
\mathcal G:=\{n_k^{\pm},\; k=1,.., k^*: |n^{\pm}_k - \eps^{-1}j\delta|
\le \eps^{-\frac 12 -\ga}\}\nn
	\end{eqnarray}
The strategy now is to fix
$\{n^{\pm}_k, k=1,..,k^*\} \in \mathcal G$
and  prove  estimates uniform  in the choice of
$\{n^{\pm}_k, k=1,..,k^*\}$, as the contribution to \eqref{5.15} of
the complement of $\mathcal G$ has negligible probability.
We have
\begin{eqnarray}
  \label{anna4.7}
&&
\max_{k=1,..,k^*} |\xi^{(\delta,-)}_{k \eps^{-2}\delta}| \le  \bar c\eps^{-1},\quad  \text{for all $\{n^{\pm}_k, k=1,..,k^*\}\in \mathcal G$}
	\end{eqnarray}
where $\bar c \eps^{-1}\ge |\xi_0| + 2 k^*\eps^{-\frac 12 - \ga} $.

\vskip.5cm

Recalling \eqref{anna4.4} for notation and that $w$ is defined in \eqref{5.11},
having fixed $\{n_k^{\pm}, \; k=1,.., k^*\} \in \mathcal G$, see \eqref{anna4.6}, with  $n^{\pm}_0\equiv 0$,
we call
\begin{eqnarray}
  \label{anna4.26}
&&  \mathcal C =  \Big\{\xi_k^0,\; k=1,.., k^* : \max_{k=1,..,k^*}\max_{x} |\mathcal A_\ell(x,\xi_k^0)- \mathcal A_\ell(x,w(\cdot|\xi_{k-1})| \le \eps^a;\nn\\&&
\hskip4cm \max_{k=1,..,k^*}  \|\xi^0_k\|_\infty\le \eps^{-a^*}\Big\}
	\end{eqnarray}
Then by  Proposition \ref{teo:3.2} and \eqref{anna4.6} after using Chebishev with the fourth power,
\begin{equation}
  \label{anna4.27}
P\Big[\{n_k^{\pm}, \xi_k^0,\; k=1,.., k^*\}\in \mathcal G \cap \mathcal C\Big] \ge 1 - c\eps^{-1-4a +2b}
= 1 - c \eps^{6/10}
	\end{equation}
The proof of \eqref{5.15.2}  continues by showing that in the set $\mathcal G \cap \mathcal C$,
$\xi_k$ (as defined in \eqref{anna4.4})
is ``close'' to $u_k$ (as in Definition \ref{definanna4.1}).  More precisely call $X_k$ and $R_k$
the integers such that
  \[
F_\eps( X_k+1;\xi^0_k) < n^+_k \le  F_\eps(X_k;\xi^0_k); \quad
F_\eps( R_k+1;u^0_k) < \eps^{-1}j\delta \le  F_\eps(R_k;u^0_k)
   \]
(see again Definition \ref{definanna4.1} for notation).  Then the
analogue of \eqref{e4anna.2} holds:

\medskip

\begin{proposition}

There are $\alpha>\beta > 1$ so that if
$\{n_k^{\pm}, \xi_k^0,\; k=1,.., k^*\}\in \mathcal G \cap \mathcal C$ then
for all $k=1,..,k^*$
\begin{equation}
  \label{anna4.28}
\max_x |\mathcal A_\ell(x,\xi^0_k)-\mathcal A_\ell(x,u^0_k)| \le \alpha^ k \eps^a,\quad  |X_k-R_k| \le \beta \alpha^ k \eps^{-1+a}
	\end{equation}

\end{proposition}

\medskip

\noindent
{\bf Proof.}  By \eqref{anna4.26}
\begin{eqnarray}
  \label{anna4.29}
&&   |\mathcal A_\ell(x,\xi^0_k)-\mathcal A_\ell(x,u^0_k)| \le \eps^a + |\mathcal A_\ell(x,u^0_k)-
A_\ell(x,w(\cdot|\xi_{k-1})|
	\end{eqnarray}

Supposing for instance that $R_{k-1}\le X_{k-1}$ we get
\begin{eqnarray}
  \label{anna4.30}
&&\hskip-1cm |w(x|\xi_{k-1}) - u^0_k(x) | =
|\sum_y p(x,y) [\xi_{k-1}(y)- u_{k-1}(y)]|\nn\\
&& \hskip1cm\le p(x,0)|n^+_{k-1}- \eps^{-1}j\delta| +
|\sum_{y < R_{k-1}}
p(x,y)[\xi^0_{k-1}(y)- u^0_{k-1}(y)]| \nn\\&& \hskip1cm
+p(x,R_{k-1}) [
\xi^0_{k-1}(R_{k-1})+ u^0_{k-1}(R_{k-1})] 
+ \sum_{R_{k-1}< y \le X_{k-1}} p(x,y) \xi^0_{k-1}(y)\nn\\
	\end{eqnarray}
By \eqref{a-e4.8}
   $$
p(x,0)|n^+_{k-1}-
\eps^{-1}j\delta| \le  \frac {c_2\eps}{\sqrt\delta} \eps^{-\frac 12-\ga}
   $$
We decompose the interval $[1,R_{k-1}-1]$ into consecutive intervals $[z_i,z'_i]$ of length $\ell$ with
the last interval which may have length $<\ell$ and get using  \eqref{a-e4.8}
  \begin{eqnarray*}
&&
  |\sum_{0<y < R_{k-1}}
p(x,y)[\xi^0_{k-1}(y)- u^0_{k-1}(y)]|
 \\&& \hskip2cm \le \sum_i \{p(x,z_i) \ell \alpha^{ k-1}\eps^a + \sum_{z_i\le y \le z'_i} |p(x,z_i)-p(x,y)| 2\eps^{-a^*}\} + \frac{c_2 \eps}{\sqrt\delta}\ell 2\eps^{-a^*}\\&&
\hskip2cm \le \frac{c_2 \eps}{\sqrt\delta} \eps^{-1} \alpha^{ k-1}\eps^a + ( \frac{c_3 \eps}{\sqrt\delta})^2 2 \eps^{-a^*}
+ \frac{c_2 \eps}{\sqrt\delta}2\eps^{-b - a^*} \le \frac{c_2}{\sqrt\delta}  \alpha^{ k-1}\eps^a
+
c \eps^{1-b-a^*}
  \end{eqnarray*}
We also have
  $$
p(x,R_{k-1}) [
\xi^0_{k-1}(R_{k-1})+ u^0_{k-1}(R_{k-1})]\le  \frac {c_2\eps}{\sqrt\delta} 2\eps^{-a^*}
   $$
By \eqref{a-e4.8}
and    \eqref{a-e4.13.0} and decomposing as before   the interval $[R_{k-1}+1,X_{k-1}]$ into consecutive intervals of length $\ell$,
  \begin{eqnarray*}
&&
\sum_{R_{k-1}< y \le X_{k-1}} p(x,y) \xi^0_{k-1}(y) \le
\sum_{R_{k-1}< y \le X_{k-1}} p(x,y) w(y|\xi_{k-2}) \\&& \hskip2cm
+ |\sum_{R_{k-1}< y \le X_{k-1}} p(x,y) [\xi^0_{k-1}(y) -w(y|\xi_{k-2})]|
\\&& \hskip1cm \le
\frac{c_2\eps}{\sqrt\delta}
 \frac{c_2  c^* }{\sqrt \delta} |X_{k-1}-R_{k-1}| +\frac{c_2 \eps}{\sqrt\delta}  |X_{k-1}-R_{k-1}| \alpha^{ k-1}\eps^a
 + ( \frac{c_3 \eps}{\sqrt\delta})^2 2 \eps^{-a^*}
+ \frac{c_2 \eps}{\sqrt\delta}2\eps^{-b - a^*}\\&&\hskip1cm
\le c \Big(\eps|X_{k-1}-R_{k-1}|+\eps^{1-b-a^*}\Big)
  \end{eqnarray*}
By collecting the above bounds and using the induction hypothesis:
  \begin{eqnarray*}
&&\hskip-1cm
|w(x|\xi_{k-1}) - u^0_k(x) | \le  \frac {c_2\eps}{\sqrt\delta} \eps^{-\frac 12-\ga}+ \frac{c_2}{\sqrt\delta}  \alpha^{ k-1}\eps^a
+
2c \eps^{1-b-a^*}+  \frac {c_2}{\sqrt\delta} 2\eps^{1-a^*}
+ c\beta \alpha^{k-1}\eps^a \\&&\hskip1cm
\le  \alpha^{ k-1}\eps^a \Big( \frac {c_2}{\sqrt\delta} \eps^{\frac 12-\ga-a}+\{\frac{c_2}{\sqrt\delta} +c\beta\}+2c \eps^{1-b-a^*-a}+
\frac {c_2}{\sqrt\delta} 2\eps^{1-a^*-a}\Big)
\\&&\hskip1cm \le \alpha^{ k-1}\eps^a \Big(\{\frac{c_2}{\sqrt\delta} +c\beta\} + \eps^{a'} C\Big)
  \end{eqnarray*}
where $a'= \min\{ \frac 12 -\ga-a,1-b-a^*-a,1-a^*-a\}>0$. 
Hence
\begin{eqnarray*}
&&   |\mathcal A_\ell(x,\xi^0_k)-\mathcal A_\ell(x,u^0_k)| \le \eps^a[1+\alpha^{k-1} \Big(\{\frac{c_2}{\sqrt\delta} +c\beta\} + \eps^{a'} C\Big)]
	\end{eqnarray*}
For $\eps$ small enough $C\eps^{a'} \le 1$, 
\begin{eqnarray}
  \label{e4anna4.29}
&&   |\mathcal A_\ell(x,\xi^0_k)-\mathcal A_\ell(x,u^0_k)| \le  \alpha^{k}\eps^a,\quad \alpha = 2+ \{\frac{c_2}{\sqrt\delta} +c\beta\}
	\end{eqnarray}
By \eqref{e4anna4.29}
\begin{equation}
  \label{anna4.28.00}
|F_\eps(x;\xi^0_k)-F_\eps(x;u^0_k)| \le (\eps^{-1}-x )  \alpha^{k}\eps^a+ 2\eps^{-b-a^*}
	\end{equation}
hence, recalling \eqref{e4anna.1},
 \begin{eqnarray*}
&&|F_\eps(R_k;u^0_k) -j\delta \eps^{-1}| \le C_+,\quad
|F_\eps(X_k;\xi^0_k)-j\delta \eps^{-1}| \le \eps^{-a^*} + \eps^{-\frac 12-\ga} \le 2\eps^{-\frac 12-\ga}
\\&& |F_\eps(X_k;u^0_k)-j\delta \eps^{-1}| \le 2\eps^{-\frac 12-\ga}+|F_\eps(X_k;u^0_k)-F_\eps(X_k;\xi^0_k)| \le 2\eps^{-\frac 12-\ga}+ \eps^{-1+a}\alpha^k+ 2\eps^{-b-a^*}\\&&
C_- |R_k-X_k|\le |F_\eps(R_k;u^0_k)-F_\eps(X_k;u^0_k)| \le C_++2\eps^{-\frac 12-\ga}+\eps^{-1+a}\alpha^k ++ 2\eps^{-b-a^*}
 	\end{eqnarray*}
which proves \eqref{anna4.28} with  $\beta= C_-^{-1}(5+C_+)$.   \qed

 \vskip.5cm

\noindent	
{\bf Proof of
Theorem \ref{teo:3.1}.}  We need to prove \eqref{5.15.2}.  By \eqref{anna4.27}
we can reduce to configurations in $\mathcal G \cap \mathcal C$ and want to prove
that in such a set
   \begin{equation}
   \label{anna5.15.0}
\lim_{\eps \to 0} \max_{x\in [0,\eps^{-1}]} \;\eps |F_\eps(x;\xi_k) - F_\eps( x;u_k)|=  0
    \end{equation}
Let us suppose for the sake of definiteness that $R_k \le X_k$.  Then for $x\le R_k$
   \begin{equation*}
   |F_\eps(x;\xi_k) - F_\eps( x;u_k)| \le \Big|\sum_{y=x}^{R_k-1} (\xi^0_k-u^0_k)\Big| + \sum_{y=R_k}^{X_k} \xi^0_k
   +u^0(R_k)
   %
%
    \end{equation*}
Calling $\bar R_k \le R_k$ the largest integer so that $\bar R_k-x$ is a multiple integer of $\ell$, we get from
 \eqref{anna4.28}:
   \begin{equation*}
 \Big|\sum_{y=x}^{R_k-1} (\xi^0_k-u^0_k)\Big| \le (\bar R_k-x) \alpha^k\eps^{a} +2 \eps^{-a^*} \eps^{-b}
 \le \alpha^k\eps^{-1+1/20} + 2\eps^{-1+1/10 - 1/100}
    \end{equation*}
Call $\bar X_k$ the smallest integer $\ge X_k$ such that    $\bar X_k- R_k$ is a multiple integer of $\ell$, then
   \begin{equation*}
 \sum_{y=R_k}^{X_k} \xi^0_k \le |\bar X_k- R_k|\Big( \frac{c_2  c^* }{\sqrt \delta} + \eps^a\Big)
 +\ell \eps^{-a^*} \le c \{\eps^{-1 +a} +\eps^{-b-a^*}\} \le 2c \eps^{-1+\frac 1{20}}
       \end{equation*}
Analogous bounds hold for $x>R_k$ and
 \eqref{anna5.15.0} then  follows.  \qed

\vskip2cm

\section{Realization of the process}
\label{quinta}

Following {\cite{FG}}   we introduce  a graphical
construction of the process.
It is also convenient to enlarge the physical space $[0, \eps^{-1} ]$
by adding two extra sites $\{-1,\eps^{-1}+1\}$ so
that configurations $\xi$ are functions on $[-1,\eps^{-1}+1]$.
We denote by $\mathcal X$ the subset of all configurations $\xi$
such that $\xi(-1)=\infty$  while $\xi(x)$ is finite
for all $x\in [0,\eps^{-1}+1]$
By default in the sequel $\xi$ denotes elements of $\mathcal X$, thus $\xi$
is determined by its values for $x\ge 0$.
Physical configurations are recovered by restricting $\xi$  to $[0,\eps^{-1}]$.
We shall often work in the sequel with labeled
particles:

\begin{definition} [Ordered configurations in the enlarged space]
\label{defin:k5.1}
We denote by $\mathcal X^{\rm ord}$ the space of ordered sequences
$\und x=(x_1,x_2,..,x_n,..)$, $x_i\ge x_{i+1}$,
with values on  $[-1,\eps^{-1}+1]$,   such that
there are finitely many entries with $x_i\ge 0$, their number is denoted by $N(\und x)$, so that
$x_i=-1$ for $i>N(\und x)$ and $x_i\ge 0$ for $i\le N(\und x)$. We also define $M(\und x)$ as
the largest integer $n$ such that $x_n=\eps^{-1}+1$.
To each $\und x$ we associate the configuration $\xi_{\und x}\in \mathcal X$
	\begin{equation}
\label{1.11111111111}
\xi_{\und x}(x)=\sum_{i\ge 1} \mathbf 1_{x_i=x}\quad \text{for all $x\in [0,\eps^{-1}+1]$}, \qquad \xi_{\und x}(-1)= \infty
	\end{equation}
Viceversa, given any $\xi\in\mathcal X$ we define $\und x_\xi$ by labeling  the particles of $\xi$  consecutively starting from the right.  Finally, given a sequence $\und y$ with finitely many entries in
$[0,\eps^{-1}+1]$, say $y_{i_1}...y_{i_k}$, its re-ordering is the sequence $\und x$ where $x_1$
is the largest element in $y_{i_1}...y_{i_k}$, $x_2$ the second largest and so on; $x_{n}=-1$
for $n\ge k+1$.

\end{definition}

\medskip

\noindent
We shall be exploiting the fact that the physically relevant quantities are
the unlabeled configurations and we are therefore free to  label
the particles as we like.

\medskip

\begin{definition} [The probability space $(\Om,P)$]
\label{defin:k5.2}  We set
\[\Om= \prod_{i\ge 0} \Om_i,\quad
P= \prod_{i\ge 0} P_i
\]
where $\Om_i = \{ \om_i= (\und t_i,\und \si_i)\}$,
$\und t_i=(t_{1;i},\ t_{2;i}, ...)$
are infinite sequences
of increasing positive ``times'' $t_{k;i}$ and
$\und \si_i=(\si_{1;i},\si_{2;i},...)$ infinite sequences of  {symmetric ``jumps''},  $\si_{k;i} = \pm 1$.
For $i\ge 1$ $P_i$ is the  product probability {law} of a Poisson
{process} of intensity $1$ for
the time sequences $\und t_i$ and of a
Bernoulli
{process} {with parameter $1/2$} for
the {jump} sequences $\und \si_i$.  $(\Om_0,P_0)$ is the   probability space introduced in Definition
\ref{defin:anna3.1}.

\end{definition}

\medskip

\noindent
{\em Graphical representation.}
For each label $i\ge 0$ we draw a vertical time axis $\mathbb R_+$ (called the $i$-th time axis)
and on each of them we put ``marks'' (with values $\pm$) as described below.
For any
element $\om_i\in \Om_i$, $i\ge 1$, we draw on the $i$-th time axis a sequence
of  arrows,  at heights  $t_{k;i}$
pointing to  right or  left  if
$\si_{k;i}= \pm 1$ respectively (the $\si_{k;i}$ are called marks).
The marks on the 0-time axis are specified by $\om_0$: they
are $+$ or $-$ crosses
which are put at the times $t_{k;0}$
with $\pm$ being the value of $\si_{k;0}$.  To each arrow we associate
a displacement operator and to each cross a creation or annihilation operator.
Roughly speaking an arrow on the $i$-th axis indicates the displacement at that time of the
$i$-th particle, provided it is in $[0,\eps^{-1}]$ before and after the displacement (otherwise the displacement is canceled).
The creation operator moves a particle from $-1$ to $0$, while the annihilation operator takes to $\eps^{-1}+1$
the rightmost particle in $[0,\eps^{-1}]$ (if such a particle exists, otherwise the operation aborts).
The precise definitions are given below:


\medskip

\begin{definition} Creation, annihilation  and displacement operators on  $\mathcal X^{\rm ord}$,
denoted respectively by $a_0^{\pm}$ and $a_i^{\pm}$, $i\ge 1$.

\begin{itemize}

\item Let $i\ge 1$.  Then
$a_i^{\pm} \und x=\und x$, $\und x \in \mathcal X^{\rm ord}$,
if $x_i=-1$ or if $x_i=\eps^{-1}+1$. If instead $x_i\in [0,\eps^{-1}]$ then
$a_i^{\pm} \und x$
is the
re-ordering (see Definition \ref{defin:k5.1}) of $\und y$
where $y_j=x_j$ for  $j\ne i$ and
$y_i = x_i\pm 1$ if $ x_i\pm 1 \in [0,\eps^{-1}]$
while $y_i=x_i$ if  $ x_i\pm 1 \notin [0,\eps^{-1}]$.

\item   $a_0^{+} \und x =:\und y^{+}$ is defined as follows:  $y^+_j= x_j$ for $j\ne
k\equiv N(\und x)+1$ and $y_k=0$, (see Definition \ref{defin:k5.1}). Thus $N(a_0^{+} \und x)= N(\und x)+1$.

\item $a_0^{-} \und x=:\und y^{-}$ is defined as follows:  $\und y^-=\und x$ if $N(\und x)= M(\und x)$
(i.e.\ no $x_i\in [0,\eps^{-1}]$), otherwise
let $m:=M(\und x)+1\le N(\und x)$, so that $x_m \in [0,\eps^{-1}]$. Then
$y^-_m=\eps^{-1}+1$; while $y^-_j=x_j$ for $j\ne m$; thus $N(a_0^{-} \und x)= N(\und x)$
and  $M(a_0^{-} \und x)= M(\und x)+1$.

\end{itemize}

\end{definition}

\medskip

\noindent
The enlarged space has been introduced to make simpler the proof of the inequalities
of the next section, but in the end what is relevant is the restriction
$\und x \cap [0,\eps^{-1}]$ of the configuration to the physical space.
To this end we shall use the following lemma:

\begin{lemma}
Let $\und x$ and $\und x'$ be such that $N(\und x)=N(\und x')$ and $M(\und x)=M(\und x')$ then
\begin{equation}
\label{d0}
N(a_i^\si\und x)=N(a_i^\si\und x');\;\; M(a_i^\si\und x)=M(a_i^\si\und x');\;\;\quad \text{for any
$i\ge 0$ and any $a_i^\si$}
\end{equation}
\begin{equation}
\label{d0.1}
N(a_i^\si\und x)=N(\und x);\;\; M(a_i^\si\und x)=M(\und x);\;\;\quad \text{for any $i\ge 1$ and any $a_i^\si$}
\end{equation}
and for any sequence $a_{i_j}^{\si_{i_j}}$, $j=1,..,n$,
\begin{equation}
\label{d0.3}
N\Big(\prod_{j=1}^n a_{i_j}^{\si_{i_j}} \und x\Big)=
N\Big(\prod_{j=1}^n (a_{i_j}^{\si_{i_j}})^{\mathbf 1_{i_j=0}} \und x\Big),\;\;\;
M\Big(\prod_{j=1}^n a_{i_j}^{\si_{i_j}} \und x\Big)=
M\Big(\prod_{j=1}^n (a_{i_j}^{\si_{i_j}})^{\mathbf 1_{i_j=0}} \und x\Big)
\end{equation}

\end{lemma}

\medskip
\noindent
{\bf Proof.}  $a_i^\si\und x$, $i\ge 1$, differs
from $\und x$ only if $x_i \in [0,\eps^{-1}]$
and in such a case it is obtained by rearranging
the particles in $\und x \cap [0,\eps^{-1}]$, hence \eqref{d0.1}.
Thus \eqref{d0} is a consequence of \eqref{d0.1} for $i\ge 1$.
When $i=0$, $a_0^+\und x$ increases $N(\cdot)$ by 1 leaving $M(\cdot)$ unchanged.
$a_0^-\und x=\und x$ if $N(\und x)=M(\und x)$  while if $N(\und x)>M(\und x)$ then  $M(a_0^-\und x)
=M(\und x)+ 1$,   $N(a_0^-\und x)=N(\und x)$. \eqref{d0.3} follows by applying repeatedly
\eqref{d0}.

  \qed

\medskip

\begin{definition}
\label{defin5.4mm}

Fix  $t>0$.  Then
with $P$ probability 1 $\und t_0\cap [0,t]$ has finitely many elements
which are all mutually
distinct. We define
\begin{equation}
\label{d1}
C_t(\om_0) = {\rm card}\, \Big\{ t_{k,0}\in \und t_0: t_{k,0} \le t,\; \si_{k;0}=+\Big\}
\end{equation}
and given $\und x\in \mathcal X^{\rm ord}$ let $n \ge C_t(\om_0)+N(\und x)$.
Thus it
is  well defined (with $P$ probability 1)
the sequence $\und t=(t_1,..,t_k)$, $0\le  t_j <t_{j+1}\le t$ of all times
$t_{k;i} \in [0,t]$, $k\ge 1$, $i=0,..,n$. We call $i_j$, $j=1,..,k$, the label of the time
axis to which $t_j$ belongs and $\si_j$ the corresponding $\pm$ mark.

\end{definition}

\medskip

\begin{definition} [The time flows] $T^0_t(\und x,\om)$ and $T_t(\und x,\om)$,
$t>0$, 
$\und x\in \mathcal X^{\rm ord}$ and $\om \in \Om$, are defined ($P$ almost surely)  as
follows. Let
$\und t$ be as in the previous definition, then  using the
same notation,
\begin{equation}
\label{d1}
    T^0_t(\und x,\om) = \prod_{i=1}^k (a_i^{\si_i})^{\mathbf 1_{i> 0}} \und x,\qquad
    T_t(\und x,\om) =  \prod_{i=1}^k a_i^{\si_i} \und x
\end{equation}
To define $T^{(\delta,\pm)}_{N\delta\eps^{-2}}(\und x,\om)$, $N$ a positive integer, we split
$\und t$ (defined as in Definition \ref{defin5.4mm} with $t\to N\delta\eps^{-2}$)
in $N$ groups: $\und t^{(1)},.., \und t^{(N)}$ where $\und t^{(h)}= \und t \cap
[(h-1)\eps^{-2}\delta,h\eps^{-2}\delta]$
(with $P$ probability 1 we may suppose that all such times are mutually distinct).  We then set
\begin{eqnarray}
\label{d2}
   &&
   T^{(\delta,-)}_{N\delta\eps^{-2}}(\und x,\om) = \prod_{h=1}^N\Big\{
 \prod_{i=1}^{k_h}(a^{\si^{(h)}_i}_i)^{\mathbf 1_{i= 0}} \prod_{i=1}^{k_h}
 (a^{\si^{(h)}_i}_i)^{\mathbf 1_{i> 0}} \Big\}\und x
\end{eqnarray}
\begin{eqnarray}
\label{d3}
   &&
   T^{(\delta,+)}_{N\delta\eps^{-2}}(\und x,\om) = \prod_{h=1}^N\Big\{
 \prod_{i=1}^{k_h}
 (a^{\si^{(h)}_i}_i)^{\mathbf 1_{i> 0}}  \prod_{i=1}^{k_h}(a^{\si^{(h)}_i}_i)^{\mathbf 1_{i= 0}}\Big\}
 \und x
\end{eqnarray}
We finally define $T^{(\delta,-)}_{t}(\und x,\om)$, $t\in
((N-1)\eps^{-2}\delta,N\eps^{-2}\delta)$ by dropping from the product in \eqref{d2} all operators
of the last group with $t_i^{(N)} >t$, $i\ge 1$, as well as all the creation-annihilation operators of $\und t^{(N)}$.
Also for $T^{(\delta,+)}_{N\delta\eps^{-2}}(\und x,\om)$, $t\in
((N-1)\eps^{-2}\delta,N\eps^{-2}\delta)$ we drop from the product in \eqref{d3}  all operators
of the last group with $t_i^{(N)} >t$, $i\ge 1$, but we retain
the creation-annihilation operators of $\und t^{(N)}$.

In other words in $T^{(\delta,+)}_{N\delta\eps^{-2}}$
the creation-annihilation operators of the $N$-th group occur all at time $(N-1)\delta\eps^{-2}$, while
in $T^{(\delta,-)}_{N\delta\eps^{-2}}$ they occur at time  $N\delta\eps^{-2}$, thus the above rule
for defining  $T^{(\delta,\pm)}_{t}(\und x,\om)$ means that we drop  all the
operators which appear at times larger than $t$.

\medskip

It is easy to see that the marginal over unlabeled configurations
of each one of the processes $\{T^0_t(\und x,\om)$, $T_t(\und x,\om)$,  $T^{(\delta,\pm)}_{t}(\und x,\om)\}$
has the law respectively of the free process $\xi^0_t$, the interacting process $\xi_t$
and the auxiliary processes $\xi_t^{(\delta,\pm)}$.  It also follows from \eqref{d0.3} that
\begin{eqnarray}
\label{d0.2}
   &&
N(T^{(\delta,\pm)}_{n\eps^{-2}\delta}(\und x,\om))= N(T_{n\eps^{-2}\delta}(\und x,\om)),\;\;\;
M(T^{(\delta,\pm)}_{n\eps^{-2}\delta}(\und x,\om))= M(T_{n\eps^{-2}\delta}(\und x,\om))
\end{eqnarray}

\end{definition}

\begin{figure}[h!]
\label{figura1}
\centering
\includegraphics[width=13cm]{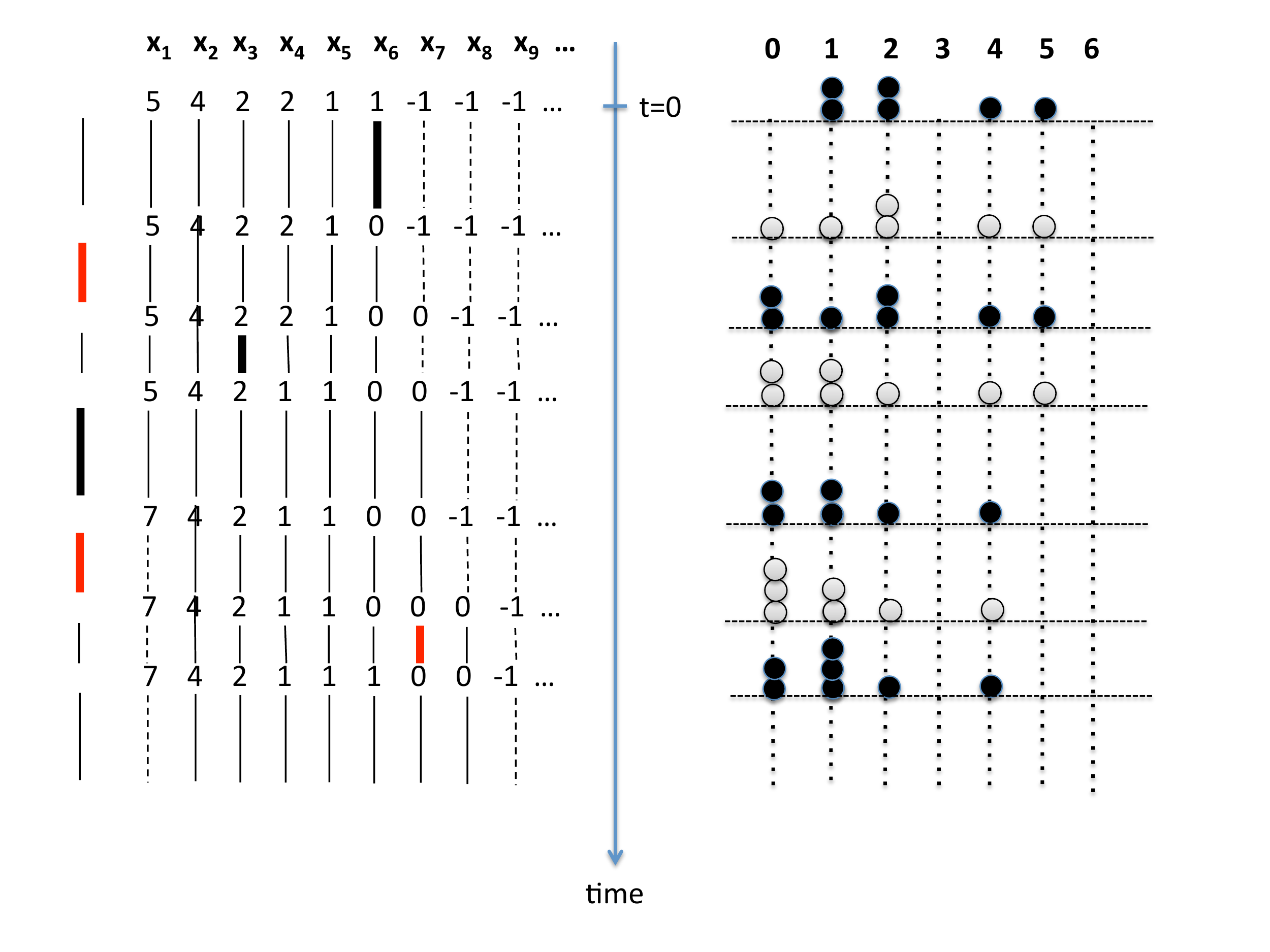}
\caption{\small{
Graphical construction of the flow $T_t$ (left panel) and of the process with generator
\eqref{generatore} (right panel) for a system of size $\eps^{-1}=6$. The legend for the
left panel is as follows: continuous vertical line denotes
the clocks of the particles involved in the dynamics; the clock of the boundaries is that
on the left;  the clock that rings first is depicted with a bold line with color red if it has associated
a jump $+1$ and color black if corresponds to a jump $-1$. After jumps particles
are re-ordered (if needed). On the right panel the motion in the physical space
$[0,6]$ is displayed.
}}
\end{figure}

%
%
%
%
%
%
%
%
%
%
%
%

\vskip2cm

\section{Mass transport inequalities}
\label{sesta}

In this section we introduce a partial order
among measures based on moving mass to the right,
we are evidently in the context of mass transport theory
from where we are borrowing the notions used in this
section.  We work
first in  the space of particle configurations
$\xi$ regarding $\xi$ as a distribution of masses and then
in the space $\mathcal U$, considering $u\in \mathcal U$ as
a mass density (which may have a Dirac delta at 0), the notions are the same
except for a change of language.

The main goal is to prove inequalities between  $\xi_t$
and the auxiliary processes $\xi_t^{(\delta,\pm)}$
(recall that the
hydrodynamic limit of the latter is known since Section \ref{quarta})
and then derive analogous inequalities for $S_t^{(\delta,\pm)}(u)$
and their limit as $\delta\to 0$.

We tacitly suppose in the sequel that the configurations $\xi$ are in $\mathcal X$
 as specified in the beginning of Section \ref{quinta}.

\begin{definition}[Partial order]
\label{coro:SumOrder}
For any $\xi,\xi'\in\mathcal X$, we say that
$\xi \le \xi'$ iff
\begin{equation}
\label{1.1.2}
F_\eps(x;\xi) \le F_\eps(x;\xi')\quad \text{for all $x\in [0,\eps^{-1}+1]$}\;.
\end{equation}


\end{definition}

\medskip
\noindent

Observe that $\xi\le \xi'$ has not the usual meaning, i.e.\ $\xi(x)\le \xi'(x)$ for all $x$ !\;
The notion of order has rather to be interpreted in the sense of ``the interfaces''
$F_\eps(x;\xi)= \sum_{y\ge x} \xi(y)$, see  Definition \ref{defin:ee8.1}
and
Figure 2 for a visual illustration.
One can easily check that the above ``$\le$'' relation
has indeed all the properties
of a partial order.
Same considerations apply to the case of continuous mass distributions as
in \eqref{e4.3.0.2.1.000} where the notion is well known and much used in mass transport theory.
\begin{figure}[h!]
\label{fiuguradue}
\centering
\vspace{-5cm}
\includegraphics[width=15cm]{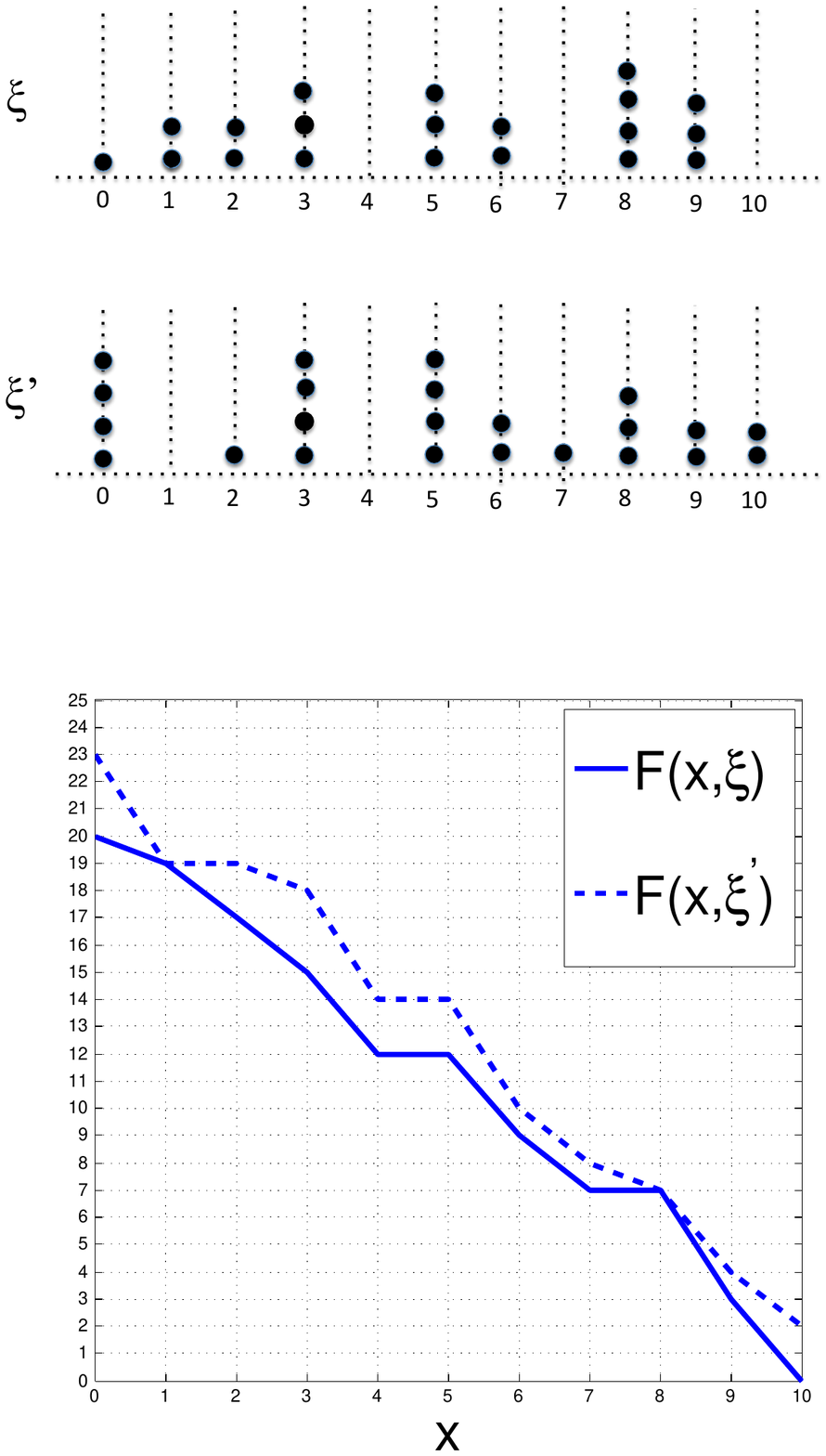}
\vspace{-1.5cm}
\caption{\small{An example of two particle configurations $(\xi,\xi')$ related by the inequality $\xi\le\xi'$ for $\eps^{-1}=10$.
Note that for the sites $x\in\{1,2,8,9\}$ one has $\xi(x) > \xi'(x)$. However the interface
of $\xi$ is below the interface of $\xi'$ for all $x\in [0,10]$.}}
\end{figure}


The equivalence with the previous statement about moving mass
to the right is established next. 
We first introduce a partial order in  $\mathcal X^{\rm ord}$ by saying that $\und x\le \und x'$ iff
$x_i\le x'_i$ for all $i$. Since there is a one-to-one correspondence between $\mathcal X$ (see Definition \ref{coro:SumOrder})  and $\mathcal  X^{\rm ord}$ this defines a priori a new order in  $\mathcal X$, but the two orders are the same as proved  in the following Proposition.

\medskip

\begin{proposition}
\label{thmk5.2}
The  conditions:   
(1)\;
$\xi \le \xi'$; (2)\;
$\und x_\xi \le \und x_{\xi'}$ (see Definition \ref{defin:k5.1})  are equivalent.
Moreover, let  $\und x=(x_1,..,x_m)$ and $\und x'=(x'_1,..,x'_n)$ be sequences with values in
$[0,\eps^{-1}+1]$
then
$\xi_{\und x} \le \xi_{\und x'}$ (see \eqref{1.11111111111}) iff
$n\ge m$ and
there is a one to one map $i_j$ from $\{1,..,m\}$ into $\{1,..,n\}$ so that
$x'_{i_j}\ge x_j$ for all $j=1,..,m$.

%

\end{proposition}

\medskip

\noindent
{\bf Proof.} Equivalence of (1) and (2). Shorthand $\und x= \und x_\xi$,  $\und x'= \und x_{\xi'}$.

Suppose  (2) holds, then
\begin{equation}
\label{k5.16}
F_\eps(x;\xi)=\sum_{i\ge 1} \mathbf 1_{x_i \ge x} \le
\sum_{i\ge 1} \mathbf 1_{x'_i \ge x}=F_\eps(x;\xi')\quad \text{for all $x\ge 0$}
\end{equation}
hence  $(2) \Rightarrow (1)$.

Suppose (1) holds and let
$\und x=(x_1,..,x_m)$ and $\und x'=(x'_1,..,x'_{n})$.  Then
$n\ge m$ because otherwise $F_\eps(0;\xi) > F_\eps(0;\xi')$.  We also have that $x_i\le x'_i$ for
$i\le m$: suppose by contradiction that $x_k >x'_k$ then $F_\eps(x_k;\xi) \ge k$
while $F_\eps(x_k;\xi') < k$, hence the contradiction.  Thus $(1) \Rightarrow (2)$.

Let  $\und x=(x_1,..,x_m)$ and $\und x'=(x'_1,..,x'_n)$ be sequences with values in
$[0,\eps^{-1}+1]$
such that
$n\ge m$ and with a one to one map $i_j$ as in the text of the proposition.
Then
\begin{equation}
\label{k5.17}
F_\eps(x;\xi_{\und x})=\sum_{j\ge 1}^m \mathbf 1_{x_j \ge x} \le \sum_{j\ge 1}^m \mathbf 1_{x'_{i_j} \ge x}
\le
F_\eps(x;\xi_{\und x'})
\end{equation}
hence $\xi_{\und x} \le \xi_{\und x'}$.
To prove the converse statement
we suppose that $\und x=(x_1,..,x_m)$ and $\und x'=(x'_1,..,x'_n)$  are such that
$\xi:=\xi_{\und x} \le \xi':=\xi_{\und x'}$.
Then $\und y:= \und x_\xi\le \und y' = \und x_{\xi'}$, and there are one to one maps
$\ell_j:\{1,..,m\}$ onto itself and $\ell'_j:\{1,..,n\}$ onto itself so that
$y_{\ell_j} = x_j$ and $x'_{\ell'_h}= y'_{h}$. Then $x_j\le x'_{i_j}$
with $i_{j} = \ell'_{\ell_j}$.

  \qed

  \medskip

As a corollary we have

\medskip

  \begin{lemma}
   \label{lemmak67}
If $\und x \le \und x'$ then
\begin{eqnarray}
\label{k5.18}
&&a^{\pm}_i \und x \le a^{\pm}_i \und x', \;\; i\ge 0;\qquad\und x \le  a^{\pm}_0\und x'\nn
\\&&\\
&& a^{+ }_0 \und x \le \und x' \;\; \text{if }\; N(\und x) < N(\und x')\qquad
a^{-}_0 \und x \le \und x' \;\; \text{if }\; M(\und x) < M(\und x')
\nn
\end{eqnarray}

   \end{lemma}

\medskip

\noindent
{\bf Proof.} The inequality $\und x \le  a^{\pm}_0\und x'$ holds trivially
because $a^{\pm}_0\und x'$ does not decrease the entries of $\und x'$.
Let us next consider the other inequalities involving  $a^+_0$. Let $k=N(\und x)+1$, then
$\und y:=a^+_0 \und x$ has $y_k=0$, while $x_k=-1$ (all the other entries are unchanged).
If $N(\und x')>N(\und x)$ then $x'_{k}\ge 0$ and the last inequality
in  \eqref{k5.18} is satisfied. If $N(\und x')=N(\und x)$ then $x'_k=-1$
but $y'_k=0$, where $\und y'= a^+_0 \und x'$, hence the first equality in
\eqref{k5.18}.

Let us next consider $a^-_0$. If $M(\und x)=N(\und x)$ then $a^-_0\und x=\und x$
and therefore  is $\le \und x' \le a_0^-\und  x'$. Let then $m=M(\und x)+1 \le N(\und x)$.
Then $\und y:= a^-_0\und x$ has $y_m = \eps^{-1}+1$. If $M(\und x')>M(\und x)$, then
$x'_m=\eps^{-1}+1$.  If instead $M(\und x')=M(\und x)$ then $x_m\le x'_m$ hence $x'_m \in [0,\eps^{-1}]$
and $\und y' = a^-_0\und x'$ has $y'_m=\eps^{-1}+1$.

Let next $\und y =a^{\pm}_i \und x$ and $\und y'= a^{\pm}_i \und x'$ with $i\ge 1$
and for the sake of definiteness let us just consider the $+$ case.
$\und y = \und x$ if $i \le M(\und x)$ and $i>N(\und x)$.  In the former case
$x'_i = \eps^{-1}+1$ is also unchanged, in the latter $x_i=-1$ and again the inequality holds
trivially.  Let us then suppose that $ M(\und x) < i\le N(\und x)$ and suppose that
this holds as well for $\und x'$ (otherwise $x'_i=\eps^{-1}+1$).
Then  $\min\{x_i + 1, \eps^{-1}\} \le \min\{x'_i + 1, \eps^{-1}\}$ hence
the desired inequality applying the last statement in Proposition
\ref{thmk5.2}.  \qed

\medskip

As already mentioned we ultimately need inequalities for the restrictions
$\und x \cap[0,\eps^{-1}]$ of the configurations to the physical space.
We shall use the following simple observation:

\medskip

\begin{lemma}
If $\und x\le \und x'$ then $N(\und x) \le N(\und x')$ and
$M(\und x) \le M(\und x')$, however
$(\und x \cap[0,\eps^{-1}])\le (\und x' \cap[0,\eps^{-1}])$
requires that $M(\und x)=M(\und x')$.  In particular if  $\und x\le \und x'$:
 \begin{equation}
\label{d17}
 (\und x \cap[0,\eps^{-1}])\le (\und x' \cap[0,\eps^{-1}])\;\; {\rm if}\;\;
 N(\und x) = N(\und x'),\,\;\; M(\und x) = M(\und x')
\end{equation}

\end{lemma}

\medskip
\begin{definition}[Stochastic order]
A process $(\xi_t)_{t\ge 0}$  is stochastically  smaller than
a process  $(\xi'_t)_{t\ge 0}$, writing in short $\xi_t \le \xi'_t$ (stochastically),
if they can be both realized on a same space where the inequality holds pointwise almost surely.

\end{definition}

\medskip

We shall prove stochastic
order by realizing the   processes on the same space $(\Om, P)$ of Definition \ref{defin:k5.2}.\

\begin{definition}
A  map $f: \mathcal X^{\rm ord} \to \mathcal X^{\rm ord}$ preserves order if
$\und x\le \und x'$ implies $f(\und x) \le f(\und x')$.

\end{definition}

\medskip
The first inequality in \eqref{k5.18} proves that all the maps $a_i^{\pm}$ preserve order
and since all the flows have been defined in terms of products of such maps:

\medskip

  \begin{theorem}[Stochastic inequalities]
  \label{Dis-}
 All the maps $T^{(\delta,\pm)}_{m\eps^{-2}\delta}(\cdot,\om)$, $T^{0}_{t}(\cdot,\om)$ and
  $T_{t}(\cdot,\om)$,
preserve order.

  \end{theorem}

  \medskip

To compare the flows   $T^{(\delta,\pm)}_t$ and $T_t$ we shall use the following lemma:

\medskip

  \begin{lemma}
  Let $i\ge 1$, then
  \begin{equation}
\label{d5}
 a_0^{\si_0} a_i^{\si_i} \und x \le
 a_i^{\si_i} a_0^{\si_0} \und x
\end{equation}

  \end{lemma}

\medskip

\noindent
{\bf Proof.} Let $\si_0=+$. Call $\und y=a_0^{+} \und x$, then by the second inequality
in \eqref{k5.18}, $\und x \le \und y$.  Since $a_i^{\si_i}$ preserves order:
$ a_i^{\si_i} \und x  \le  a_i^{\si_i} \und y$ and since $N(\und y)= N(\und x)+1$ we have
\eqref{d5} (having used the third inequality
in \eqref{k5.18}).  Let $\si_0=-$. Call $\und y=a_0^{-} \und x$, then by the second inequality
in \eqref{k5.18}, $\und x \le \und y$.  Since $a_i^{\si_i}$ preserves order:
$ a_i^{\si_i} \und x  \le  a_i^{\si_i} \und y$ and since $M(\und y)= M(\und x)+1$ we have again
\eqref{d5} (having used the fourth inequality
in \eqref{k5.18}).   \qed

\medskip

  \begin{corollary}
Let $\{(i_j,\si_j)\}$ a sequence of  $n> 1$ pairs with $i_j \ge 0$, $\si_j\in \{+,-\}$.
An exchange  at $(h,h+1)$, $h+1\le n$, is the new sequence $\{(i'_j,\si'_j)\}$
where $(i'_j,\si'_j)=(i_j,\si_j)$ for $j\ne h,h+1$ and  $(i'_h,\si'_h)=(i_{h+1},\si_{h+1})$,
$(i'_{h+1},\si'_{h+1})=(i_h,\si_h)$.  We then say that an  exchange at $(h,h+1)$
is  ``allowed''
if $i_h=0$ and $i_{h+1}>0$.

Then if $\pi$ is a permutation obtained by applying repeatedly allowed exchanges starting from
 $\{(i_j,\si_j)\}$ so that the final sequence is $\{(i_{\pi(j)},\si_{\pi(j)})\}$
  \begin{equation}
\label{d6}
 \prod_{j=1}^n a_{i_j}^{\si_j}\und x \le
  \prod_{j=1}^n a_{i_{\pi(j)}}^{\si_{\pi(j)}}\und x
\end{equation}

  \end{corollary}

\vskip.5cm

Call $\{(a_j,\si_j)\}$ the sequence associated to
$T^{(\delta,-)}_{m\eps^{-2}\delta}(\und x,\om)$ and
$\{(a'_j,\si'_j)\}$ the one associated to $T_{m\eps^{-2}\delta}^{(\delta',-)}(\und x,\om)$, $\delta=k\delta'$: then the
latter is obtained by repeated allowed exchanges from the former, hence
 	\begin{equation*}
T_{m\eps^{-2}\delta}^{(\delta,-)}(\und x,\om) \le T_{m\eps^{-2}\delta}^{(\delta',-)}(\und x,\om)
 	\end{equation*}
Also the sequence $\{(a''_j,\si''_j)\}$ associated to $T_{m\eps^{-2}\delta}(\und x,\om)$
is obtained by repeated allowed exchanges from $\{(a'_j,\si'_j)\}$, hence
 	\begin{equation*}
T_{m\eps^{-2}\delta}^{(\delta',-)}(\und x,\om) \le
T_{m\eps^{-2}\delta}(\und x,\om)
 	\end{equation*}
The sequence $\{(a'''_j,\si'''_j)\}$ associated to $T_{m\eps^{-2}\delta}^{(\delta',+)}(\und x,\om)$
is obtained by repeated allowed exchanges from $\{(a''_j,\si''_j)\}$, hence
 	\begin{equation*}
T_{m\eps^{-2}\delta} \le
T_{m\eps^{-2}\delta}^{(\delta',+)}(\und x,\om)
 	\end{equation*}
Finally  the sequence $\{(a^*_j,\si^*_j)\}$ associated to $T_{m\eps^{-2}\delta}^{(\delta,+)}(\und x,\om)$
is obtained by repeated allowed exchanges from $\{(a'''_j,\si'''_j)\}$, hence
 	\begin{equation*}
T_{m\eps^{-2}\delta}^{(\delta',+)}(\und x,\om) \le
T_{m\eps^{-2}\delta}^{(\delta,+)}(\und x,\om)
 	\end{equation*}
We have thus proved:

\medskip

  \begin{theorem}[Stochastic inequalities]
  \label{Dis-com}
  Denoting by $\xi_{m\eps^{-2}\delta}^{(\delta,\pm)}$ and $\xi_{t}^{(\delta,\pm)}$
  the configurations $\xi_{T_{m\eps^{-2}\delta}^{(\delta,-)}(\und x,\om)}$
  and $\xi_{T_{t}(\und x,\om)}$ restricted to $x\in [0,\eps^{-1}]$ we have
  for any $\delta = k\delta'$, $k$ a positive integer,
 	\begin{equation}
	\label{ee8.4bis}
\xi_{m\eps^{-2}\delta}^{(\delta,-)} \le \xi_{m\eps^{-2}\delta}^{(\delta',-)}
\le \xi_{m\eps^{-2}\delta}\le
\xi_{m\eps^{-2}\delta}^{(\delta',+)} \le
\xi_{m\eps^{-2}\delta}^{(\delta,+)}
 	\end{equation}

  \end{theorem}

  \medskip

  \noindent
{\bf Proof.}  We have already proved the inequality for the configurations
on $[-1,\eps^{-1}+1]$, thus the proof of \eqref{ee8.4bis} follows from
\eqref{d17} and \eqref{d0.2}.  \qed

  \vskip.5cm

%
%

The theorem has its continuum analogue which can be proved directly, see Section 4 of \cite{CDGP3},
but it can also be deduced from Theorem \ref{Dis-com}, as we shall see.

   \medskip

\begin{theorem}[Macroscopic inequalities]
\label{thmee8.1}
Let $u\in L^\infty([0,1],\mathbb R_+)$, $F(0;u)>0$.  Let $\delta: j\delta < F(0;u)$ and
$\delta'$ such that $\delta = k\delta'$ with $k$ a positive integer.  Then
 	\begin{equation}
	\label{ee8.4}
S_{m \delta}^{(\delta,-)}(u) \le S_{m \delta}^{(\delta',-)}
\le
S_{m \delta}^{(\delta',+)} \le
S_{m \delta}^{(\delta,+)}
 	\end{equation}
Moreover  the maps $K^{(\delta)}$, ${G^{\rm neum}_t}*$  and $S_{t}^{(\delta,\pm)}$
on $\mathcal U_\delta$,  see \eqref{e4.3.0.0}, preserve  order.

	\end{theorem}

\medskip

\noindent
{\bf Proof.}  \eqref{ee8.4} follows from \eqref{ee8.4bis} and \eqref{5.15}.
Proof that $K^{(\delta)}u\le K^{(\delta)}v$, $u,v\in \mathcal U_\delta$. We have
      \begin{equation*}
K^{(\delta)}u -K^{(\delta)} v= (c_u-c_v)D_0 +  (\rho_u-\rho_v) \mathbf 1_{r\le R_\delta(u)}
-\rho_v \mathbf 1_{ R_\delta(u) < r \le R_\delta(v)}
   \end{equation*}
where $R_\delta(w): F(R_{\delta}(w);w)=j\delta$.   Hence
      \begin{equation*}
F(r;K^{(\delta)} u)-F(r; K^{(\delta)} v)= \Big(F(r;u)-F(r;v)\Big) \mathbf 1_{r \le
R_\delta(u)} - \mathbf 1_{r >
R_\delta(u)} \int_r^{R_\delta(v)}
\rho_v(r') 
   \end{equation*}
which is therefore $\le 0$.

The property that ${G^{\rm neum}_t}*$
preserves the order is inherited from the same property for
the independent flow $T^0_t$. As a consequence of the two previous statements we have that also $S^{(\delta,\pm)}_t$ preserves the order (see the definition in \eqref{e4.3.0.0}).  \qed

\vskip2cm

\section{Regularity properties of the  barriers}
\label{settima}

In this section we shall prove some regularity properties of
the  barriers $S_{t}^{(\delta,\pm)}(u)$,
$u\in L^\infty([0,1],\mathbb R_+)$, $F(0;u) >j\delta$ (the barriers are
defined in Definition \ref{defin:e2.3}).

By the smoothness of $G_{t}^{\rm neum}(r,r')$, $t>0$, it is easy
to prove that for any $n>0$,
$S^{(\delta,+)}_{n\delta}(u)$ is in $C^\infty$ while
$S^{(\delta,-)}_{n\delta}(u)$
is equal to $j\delta D_0$ plus a function which is
$C^\infty$ in the interior of its support.
Such a smoothness however, being inherited from
$G^{\rm neum}_\delta$, depends on $\delta$, while we want
properties which hold uniformly as $\delta\to 0$.

The properties of the Green functions that we use
in this section are:
 	\begin{equation}
	\label{e6.1aaa.00}
G_{t}^{\rm neum}(r,r') =G_{t}^{\rm neum}(r',r)  \le \frac{c(1+\sqrt{t})}{\sqrt{t}}, \quad
|\frac{d}{dr} G_{t}^{\rm neum}(r,r')| \le \frac{c}{t} 
	\end{equation}
 	\begin{equation}
	\label{e6.1aaa.01}
\int dr'G^{\rm neum}_{t} (r,r') = 1 
	\end{equation}
 	\begin{equation}
	\label{e6.1aaa.01.1}
\int_{|r'-r|>X} dr'G^{\rm neum}_{t} (r,r') \le \sqrt 2 e^{- X^2/(4t)},\qquad \forall X>0
	\end{equation}
\eqref{e6.1aaa.01.1} is proved by writing
 	\begin{equation*}
\int_{|r'-r|>X} dr'G^{\rm neum}_{t} (r,r') \le \int_{|r'-r|>X} dr'G_{t} (r,r'),\quad G_t(r,r') = \frac{e^{-(r-r')^2/(2t)}}{\sqrt{2\pi t}}
	\end{equation*}
and then bounding
 	\begin{equation*}
 \int_{|r'-r|>X} dr'G_{t} (r,r')\le e^{-X^2/(4t)} \sqrt 2\int dr \frac{e^{-r^2/(4t)}}{\sqrt{4\pi t}}
	\end{equation*}

Such bounds are verified also by the Green function for the Neumann problem in $[0,\ell]$
for any $\ell >0$ and $\ell=\infty$ as well, so that the analysis in this section extends
to all such cases.  Observe that if $\ell$ is finite and positive the bound on the derivative is
much better:
 	\begin{equation*}
\Big|\frac{d}{dr} G_{t}^{\rm neum}(r,0)\Big| \le \frac{ce^{-bt}}{t},\quad b>0, c>0 
	\end{equation*}
but we shall only use \eqref{e6.1aaa.00}, \eqref{e6.1aaa.01} and \eqref{e6.1aaa.01.1}
to have what follows valid also in the spatial domain  $[0,\infty)$.
%


The main results in this section are:

\vskip.5cm

\begin{theorem}[Space and time equicontinuity]
\label{thmee7.1}
Let $u\in L^\infty([0,1],\mathbb R_+)$, $F(0;u) > 0$.  Then

\begin{itemize}

\item   $F(0;S_{t}^{(\delta,\pm)}(u))= F(0;u)$ for all  $\delta>0$  such that
$F(0;u) > j\delta$ and all $t=n\delta$, $n\in \mathbb N$.

\item There is a constant $c$ so that  for any $\delta>0$: $F(0;u) > j\delta$
 \begin{equation}
	\label{e6.1}
\|S^{(\delta,+)}_{t}(u)\|_{\infty} \le
c\left\{
\begin{array}{ll}
j+ \|u\|_{\infty} & \text{for all $t \in \delta \N, \: t \le 1$}\\
j+ F(0;u)  & \text{for all  $t \in \delta \N, \: t>1$}
\end{array}
\right.
\end{equation}
Same bounds hold for $\{S^{(\delta,-)}_{t}(u) - jD_0\}$.

\item  Given any time $\si >0$ the following holds.
For any $\zeta>0$ there are $\tau_\zeta>0$
and $d_\zeta>0$ so that for any $\delta\in (0,\si)$: $F(0;u) > j\delta$,
for any $t\ge \si$ in $\delta\mathbb N$, for any $t' \in \delta\mathbb N$, $t'\in (t,t+\tau_\zeta)$ and for any $r$ and $r'$ such that $|r-r'|<d_\zeta$,
 	\begin{equation}
	\label{e6.1aaa.1}
| S^{(\delta,+)}_{t}(u)(r)-S^{(\delta,+)}_{t}(u)(r')| <\zeta,
\quad| S^{(\delta,+)}_{t'}(u)(r)-S^{(\delta,+)}_{t}(u)(r)| <\zeta
	\end{equation}

\item For all  $\delta>0$ such that
$F(0;u) > j\delta$ and all
$t>0$ in $\delta \mathbb N$  
 	\begin{equation}
	\label{e6.1aaa.2}
F\Big(0;|S_t^{(\delta,+)}(u)-S_t^{(\delta,-)}(u)|\Big)\le  4j \delta
	\end{equation}

\end{itemize}


\end{theorem}

\medskip

\noindent
{\bf Proof.}

$\bullet$\; $F(0;S_{t}^{(\delta,\pm)}(u))= F(0;u)$  because  by \eqref{e6.1aaa.01}  $G_\delta^{\rm neum}$
preserves the mass, as well as $K^{(\delta)}$, by its very definition, see
\eqref{e4.3.0.1}).
%

$\bullet$\;  Proof of \eqref{e6.1}.
Let $t=n\delta $, $n$ a positive integer, then
  \[
  S^{(\delta,+)}_{t}(u) (r) \le
  \int dr'G_{\delta}^{\rm neum}(r,r') S^{(\delta,+)}_{t-\delta}(u)(r') + j \delta
G_{\delta}^{\rm neum}(r,0)
  \]
The inequality is because we are not taking into account the ``loss part''
in the action of $K^{(\delta)}$.  Iterating we get for $s=m\delta$, $m< n$ a non negative integer,
	\begin{equation}
	\label{e6.2}
S^{(\delta,+)}_{t} (u)(r) \le \int dr'G_{t-s}^{\rm neum}(r,r') S^{(\delta,+)}_{s}(u)(r') + j \delta \sum_{k=1}^{n-m}
G_{k\delta}^{\rm neum}(r,0)
	\end{equation}
Let $n_\delta$ be the smallest integer such that
$\delta n_\delta\ge 1$ and suppose that in \eqref{e6.2}
$t < \delta n_\delta$ and  $s=0$. By \eqref{e6.1aaa.01} the integral in \eqref{e6.2}
is  bounded by $\|u\|_\infty$ 
whereas by \eqref{e6.1aaa.00} the sum  is bounded by $c''j\sqrt{n \delta}\le c''j$.
Thus  \eqref{e6.1} is proved  for $t \le 1$.

Let us next take $t=\delta n_\delta$ and $s=0$ in \eqref{e6.2}.
Then using \eqref{e6.1aaa.00} we bound the integral in
\eqref{e6.2} by $c'F(0;u) (\delta n_\delta)^{-1/2}\le c'F(0;u)$.
As before the last term in \eqref{e6.2} is
bounded by $c''j\sqrt{\delta n_\delta}\le 2c'' j$ so that (we may suppose $c'<2c''$)
\[
\|S^{(\delta,+)}_{ n_\delta \delta}(u)\|_\infty \le 2c'' (F(0;u) +j)
\]
By the same argument for any integer $k \ge 1$
   \begin{equation}
   \label{666}
\|S^{(\delta,+)}_{ k n_\delta \delta}(u)\|_\infty
\le 2c'' \{F\Big(0;S^{(\delta,+)}_{ (k-1)\delta n_\delta }(u)\Big)  +j\}=2c'' (F(0;u) +j)
   \end{equation}
the last equality because we have already proved that mass is conserved.
Thus \eqref{e6.1} is proved for $t \in (\delta n_\delta) \mathbb N$.
Let now  $m=k n_\delta$ and  $k n_\delta<n \le (k+1) n_\delta$
$k$ a positive integer.  The last term in \eqref{e6.2}
is bounded again by $2c''j$, whereas the integral
is smaller than  $\|S^{(\delta,+)}_{k n_\delta \delta}(u)\|_\infty$.
Thus \eqref{e6.1} follows from  \eqref{666} when $t\ge 1$.

We next prove the analogue of \eqref{e6.1} for
 	\begin{equation}
	\label{e6.1aaa}
 \rho^{(\delta,-)}_{t}:=
S^{(\delta,-)}_{t} (u)- j\delta D_0, 
\qquad t >0 \in \delta \mathbb N
	\end{equation}
Let  $t=n\delta$, $s=m\delta$, $n >m$  in $\mathbb N$, just as before.
Recalling the definition \eqref{e6.1aaa}, we have
  	\begin{equation*}
\rho^{(\delta,-)}_{n\delta}=K^{(\delta)} [G^{\rm neum}_\delta*
S^{(\delta,-)}_{(n-1)\delta} (\rho_0)]- j\delta D_0
= \mathbf 1_{ [0,R]}\,\,\,
G^{\rm neum}_\delta*
\big[\rho^{(\delta,-)}_{(n-1)\delta} +j\delta D_0\big]
	\end{equation*}
where $\mathbf 1_{[0,R]}$ is the characteristic function of the set $[0,R]$
and $R$ is such that
  	\begin{equation*}
\int_R^1
G^{\rm neum}_\delta*
\big[\rho^{(\delta,-)}_{(n-1)\delta} +j\delta D_0\big](r) = j\delta
	\end{equation*}
Then
 			\begin{eqnarray}
	\nn
&& \hskip-1cm\rho^{(\delta,-)}_{n\delta}(r)= \mathbf 1_{r\le R}\Big(j\delta G_\delta^{\rm neum}(r,0)+
G^{\rm neum}_\delta*
\rho^{(\delta,-)}_{(n-1)\delta}(r)\Big)\nn\\
&&\le j\delta G_\delta^{\rm neum}(r,0)+
G^{\rm neum}_\delta*
\rho^{(\delta,-)}_{(n-1)\delta}(r)
\label{merda}
	\end{eqnarray}
After iterating \eqref{merda} we get
 	\begin{equation}
	\label{e6.4a}
 \rho^{(\delta,-)}_t (r) \le  j\delta  \sum_{k=1}^{n-m} G_{k\delta}^{\rm neum} (r,0) +
 \int dr'G_{t-s}^{\rm neum}(r,r') \rho^{(\delta,-)}_s(r')
%
%
%
	\end{equation}
which has the same structure as \eqref{e6.2}.  The analysis after
\eqref{e6.2} extends to the present case and yields the proof of
\eqref{e6.1} for $\rho^{(\delta,-)}_{t}$.

\medskip

The proof of \eqref{e6.1aaa.1}
and \eqref{e6.1aaa.2} will be given after the following lemma.

\medskip

\begin{lemma}
\label{lemmae6.24}
There is a constant
$c$ so that the following holds. For all  $\delta>0$ such that
$F(0;u) > j\delta$ and for all  $0\le s< t$, $s,t \in\delta \mathbb N$, $t-s\le 1$,
we  write
 	\begin{equation}
	\label{e6.333333}
w_{s,t}^{(\delta,+)}(r) := \int dr'G_{t-s}^{\rm neum}(r,r') S^{(\delta,+)}_{s}(u)(r'),\;\;\;
v_{s,t}^{(\delta,+)}:=  S^{(\delta,+)}_{t}(u) - w_{s,t}^{(\delta,+)}
	\end{equation}
Then
%
%
%
%
%
%
 	\begin{equation}
	\label{e6.3}
\sup_{r,r' \in [0,1]}| w_{s,t}^{(\delta,+)}(r) - w_{s,t}^{(\delta,+)}(r')| \le c  \|u\|_\infty  \frac{|r-r'|}{ t-s }
	\end{equation}
 	\begin{equation}
	\label{e6.3.0}
F(0;|v_{s,t}^{(\delta,+)}|) \le 2j (t-s),\quad  
\|v_{s,t}^{(\delta,+)}\|_\infty \le c  j \sqrt{t-s}
	\end{equation}

\end{lemma}

\medskip

\noindent
{\bf Proof.}
By \eqref{e6.1} and the second inequality in \eqref{e6.1aaa.00}
we get
	\begin{equation*}
|w_{s,t}^{(\delta,+)}(r)-w_{s,t}^{(\delta,+)}(r')|\le \|S^{(\delta,+)}_{s}(u)\|_{\infty}\int  |G^{\rm neum}_{t-s}(r,z)-G^{\rm neum}_{t-s}(r',z)| \, dz\le c \frac{ |r'-r''| }{  t-s} \,\|u\|_\infty
	\end{equation*}
which proves \eqref{e6.3}.

We already have an upper bound for $S^{(\delta,+)}_{t} (u)(r)$ as given by \eqref{e6.2}
and want to find a
lower bound. We first define for any $ \tau \in \delta \mathbb N$
\begin{equation}
\label{a1}
v^{(\delta)}_\tau (r) = \mathbf 1_{r\ge R}S^{(\delta,+)}_\tau(u) (r),\quad R:\;\;
\int v^{(\delta)}_\tau (r) = j\delta
\end{equation}
By \eqref{e6.1}
\begin{equation}
\label{a2}
\|v^{(\delta)}_\tau \|_\infty \le C, \quad C=c\Big(j+ \|u\|_{\infty}\Big)
\end{equation}
By neglecting the contribution of $jD_0$ we get:
 	\begin{equation*}
S^{(\delta,+)}_{t}(u) \ge
  G^{\rm neum}_{\delta}* \Big(S^{(\delta,+)}_{t-\delta}(u) - v^{(\delta)}_{t-\delta}\Big)
	\end{equation*}
and by iteration:
 	\begin{equation*}
 S^{(\delta,+)}_{t}(u)  \ge    G^{\rm neum}_{t-s}* S^{(\delta,+)}_s(u)
 - \sum_{k= m} ^{ n-1}  G^{\rm neum}
 _{ (n- k)\delta} * v^{(\delta)}_{k\delta}
	\end{equation*}
Combining the upper and the lower bound and recalling \eqref{e6.333333}
 	\begin{equation}
	\label{e6.4aa.1}
|v_{s,t}^{(\delta,+)}|\equiv | S^{(\delta,+)}_{t}(u)  -  G^{\rm neum}_{t-s}* S^{(\delta,+)}_s(u) | \le \sum_{k= m} ^{ n-1}  G^{\rm neum}
 _{ (n- k))\delta} * v^{(\delta)}_{k\delta}+  j\delta  \sum_{k=1}^{ n-m}  G_{ k\delta}^{\rm neum} (r,0)
	\end{equation}
By \eqref{a1} and
\eqref{e6.1aaa.00} 
	\begin{eqnarray*}
\| \sum_{k=m+1}^{n}  G^{\rm neum}_{ (n- k)\delta} * v^{(\delta)}_{k\delta} \|_\infty \le cj\sqrt \delta \sqrt{n-m}= c
 j\sqrt{t-s}
 	\end{eqnarray*}
and by \eqref{e6.1aaa.00}
	\begin{eqnarray*}
\|j\delta  \sum_{k=1}^{ n-m}  G_{ k\delta}^{\rm neum} (r,0) \|_\infty \le cj\sqrt \delta \sqrt{n-m}= c
 j\sqrt{t-s}
 	\end{eqnarray*}
so that $\|v_{s,t}^{(\delta,+)}\|_\infty \le c  j \sqrt{t-s}$
and the second inequality in \eqref{e6.3.0} is proved.  To prove the first one
we use \eqref{e6.4aa.1},
\eqref{e6.333333} and  \eqref{e6.1aaa.01} to write
\begin{eqnarray*}
F(0;|v^{(\delta,+)}_{s,t}|)\le j\delta (t-s)+F\Big(0;\sum_{k=m}^{n-1}  G^{\rm neum}_{ (n- k)\delta} * v^{(\delta)}_{k\delta} \Big) \le 2j\delta(t-s)
 	\end{eqnarray*}
which concludes the proof of \eqref{e6.3.0}. 

\qed

%
%
%
%

\vskip.5cm

\noindent
We resume the proof of Theorem \ref{thmee7.1} by proving:

\medskip

$\bullet$\;\; Proof of the first inequality in \eqref{e6.1aaa.1} (space equicontinuity).
Recalling that $\delta < \si$ we
may suppose (with no loss of generality) that
  \begin{equation}
	\label{a0}
 \zeta < 2c' \sqrt {\si-\delta}, \quad
 c':=c(j + \|u\|_\infty)
  \end{equation}
with $c$ the constant in \eqref{e6.3}--\eqref{e6.3.0}.  Then, given any
such $\zeta>0$, we must find $d_\zeta>0$ so that
 	\begin{equation}
	\label{e6.3.0.1.1}
\sup_{|r-r'|< d_\zeta}\;| S^{(\delta,+)}_{t}(u)(r)-S^{(\delta,+)}_{t}(u)(r')| <\zeta,\quad  t\in \delta\mathbb N,
\; t\ge \si
	\end{equation}
By \eqref{e6.3} and \eqref{e6.3.0}
 	\begin{equation}
	\label{e6.3.0.1.2}
| S^{(\delta,+)}_{t}(u)(r)- S^{(\delta,+)}_{t}(u)(r')| \le
c' \frac{|r-r'|}{ {t-s}} +  c'\sqrt{t-s}
	\end{equation}
We shall prove \eqref{e6.3.0.1.1} with
 	\begin{equation}
	\label{e6.3.0.1.2.1}
d_\zeta < \zeta^3\min\Big\{ \frac {1}{4c'(2c')^2}\;;\;\frac{1}{c''(2c')^2} \Big\}
	\end{equation}
where $c''$ is a constant which will be specified later.

We first consider the case when $(2c')^2 \delta <\zeta^2$.  We then
choose $s<t$ as the smallest time in $\delta\mathbb N$ such that $2c'\sqrt{t-s} <\zeta$. Since $t-s= k\delta$,
for $s$ to exist it must be that $(2c')^2 \delta <\zeta^2$ which is indeed the case
presently considered.  On the other hand by \eqref{a0}, $s\ge \delta$.
Then, by the minimality of $s$, $2c'\sqrt{t-s+\delta} \ge\zeta$ so that
\[
2(t-s) \ge t-s+\delta \ge \frac{\zeta^2}{(2c')^2}
\]
By choosing  $d_\zeta$ as in \eqref{e6.3.0.1.2.1} the first term on the right hand side of \eqref{e6.3.0.1.2} is bounded by
\[
c' \frac{2(2c')^2}{\zeta^2}d_\zeta < \frac \zeta 2
\]
hence $|S^{(\delta,+)}_{t}(u)(r)- S^{(\delta,+)}_{t}(u)(r')| <\zeta$.

It remains to consider the case when  $(2c')^2 \delta \ge \zeta^2$.
Observe that
 	\begin{equation}
	\label{e6.3.0.1.3}
 S^{(\delta,+)}_{t}(u)= G^{\rm neum}_\delta*(j\delta D_0 + v),\quad v = \mathbf 1_{r\le R} S^{(\delta,+)}_{t-\delta}(u)
	\end{equation}
where $R$ is such that $\int_{R}^{1}S^{(\delta,+)}_{t-\delta}(u) = j\delta$. Hence
by \eqref{e6.1aaa.00} the space-derivative
of $ S^{(\delta,+)}_{t}(u)(r)$  is bounded by
\[
\frac{c}{\delta}\Big(j\delta + F(0;S^{(\delta,+)}_{t-\delta}(u))\Big)=: \frac{c''}{ \delta}
\]
with $c''= c(j\delta+F(0;u))$, having used that $F(0;S^{(\delta,+)}_{s}(u))=F(0;u)$.

By \eqref{e6.3.0.1.2.1} we then get
 	\begin{equation}
	\label{e6.3.0.1.4}
 | S^{(\delta,+)}_{t}(u)(r)- S^{(\delta,+)}_{t}(u)(r')| \le c''\delta^{-1}|r-r'| \le c''(\frac{\zeta^2}{(2c')^2 })^{-1}
 d_\zeta <\zeta
	\end{equation}

\medskip

$\bullet$\;\; Proof of the second inequality in \eqref{e6.1aaa.1} (time equicontinuity).
Let $t'>t\ge \si$, $t'-t\le 1$. Then 
by \eqref{e6.4aa.1} with $t \to t'$ and $s\to t$,
 	\begin{eqnarray*}
| S^{(\delta,+)}_{t'}(u)  -  G^{\rm neum}_{t'-t}*S^{(\delta,+)}_{t}(u)| &\le &
 \sum_{k= m} ^{ n-1}  G^{\rm neum}
 _{ (n- k))\delta} * v^{(\delta,+)}_{k\delta}
+ j\delta  \sum_{k=1}^{n-m}  G_{ (n- k)\delta}^{\rm neum} (r,0)\\
 &\le & cj\sqrt{t'-t}
	\end{eqnarray*}
Hence calling $\zeta'=\zeta/4$  and with $C\ge \|S^{(\delta,+)}_t(u)\|_\infty$ (see \eqref{e6.1}),
 	\begin{eqnarray*}
| S^{(\delta,+)}_{t'}(u)(r)  - S^{(\delta,+)}_{t}(u)(r) | \le \int_{r':|r-r'|\ge d_{\zeta'}} C G^{\rm neum}_{t'-t}(r,r') \,dr' +  \zeta' + cj\sqrt{t'-t}
	\end{eqnarray*}
We choose $\tau_\zeta = a \zeta^8$, $a$ a positive constant whose value will be specified later.  If
$\delta>\tau_\zeta$ there is no $t': t<t'<t+\tau_\zeta$ and the second inequality in \eqref{e6.1aaa.1}
is automatically satisfied.  Let then $\delta\le\tau_\zeta$.
We choose $a$ so that $cj\sqrt{a\zeta^4}<\zeta'$. By the decay properties of the Green function,
see \eqref{e6.1aaa.01.1},
 	\begin{eqnarray*}
\int_{r':|r-r'|\ge d_{\zeta'}}  G^{\rm neum}_{t'-t}(r,r') \,dr' \le \sqrt 2 e^{-c d_{\zeta'}^2/(4\tau_\zeta)}
	\end{eqnarray*}
Since $d_\zeta = c \zeta^3$ (see the  proof of space continuity) for $a$ small enough the
above integral
is $<\zeta'$ as well.

\vskip.5cm

We shall resume the proof of Theorem \ref{thmee7.1} after the following lemma:

\noindent

\medskip

\begin{lemma}
\label{lemma9}
Let $u$ and $v$ be both in $\mathcal U_\delta$, see \eqref{e4.3.0.0}, then
\begin{equation}
\label{e6.17}
F\Big(0;|K^{(\delta)}u-K^{(\delta)} v|\Big) \le F\Big(0;|u-v|\Big),\quad
F\Big(0;|K^{(\delta)}u-u|\Big) \le 2j\delta
\end{equation}

\end{lemma}

\medskip

\noindent
{\bf Proof.}  Supposing $R_\delta(u)>R_\delta(v)$, see \eqref{e4.3.0.1},
\begin{eqnarray*}
F\Big(0;|K^{(\delta)}u-K^{(\delta)} v|\Big)  &=& \int_0^{R_\delta(v)} |u-v| + \int_{R_\delta(v)}^{R_\delta(u)} u
\\&=& F\Big(0;|u-v|\Big) + \int_{R_\delta(v)}^{R_\delta(u)} (u-|u-v|) -  \int_{R_\delta(u)}^{\infty} |u-v|
\end{eqnarray*}
We have
\[
\int_{R_\delta(u)}^{\infty} |u-v| \ge |\int_{R_\delta(u)}^{\infty} (u-v)| = j\delta -  \int_{R_\delta(u)}^{\infty} v =  \int_{R_\delta(v)}^{R_\delta(u)} v
\]
so that
\[
F\Big(0;|K^{(\delta)}u-K^{(\delta)} v|\Big) \le F\Big(0;|u-v|\Big) - \int_{R_\delta(v)}^{R_\delta(u)}(v-u+ |u-v|)\le F\Big(0;|u-v|\Big)
\]

The second inequality in \eqref{e6.17}  follows because
\[
K^{(\delta)}u-u = j \delta D_0 - \mathbf 1_{r>R_\delta(u)} u
\]

\qed

\medskip
$\bullet$ \;
 Proof of \eqref{e6.1aaa.2}.
 The proof is actually a corollary of Lemma \ref{lemma9}  and the maximum principle
 \[
F \Big(0;|G_t^{\rm neum} *u-G_t^{\rm neum} *v\Big) \le F\Big(0;|u-v|\Big)
 \]
Shorthand $G$ for the operator $G_\delta^{\rm neum} *$ and
 \[
\phi:= K^{(\delta)}G \cdots  K^{(\delta)}G u,\qquad \psi:=G  K^{(\delta)}\cdots G K^{(\delta)} u
 \]
so that we need to bound the total variation of $\phi-\psi$. Call
 \[
 v = K^{(\delta)}u,\qquad v_n= G K^{(\delta)}\cdots G v,\qquad u_n=G K^{(\delta)}\cdots G u
 \]
Thus $u_n$ and $v_n$ are obtained by applying $G (K^{(\delta)}G)^{n-1}$ to $u$ and respectively $v$.
Since $G (K^{(\delta)}G)^{n-1}$  is a contraction we get, using \eqref{e6.17},
 \begin{eqnarray*}
F(0; |\psi-\phi|) &\le& F(0;|K^{(\delta)}u_n - v_n|) \le F(0;|K^{(\delta)}u_n - u_n|) +
F(0;|v_n-u_n|)\\
  &\le& 2j\delta + |v_n-u_n|_1 \le 2j\delta +|u-v|_1 \le 4j\delta
 \end{eqnarray*}

\medskip

\noindent
The proof of   Theorem \ref{thmee7.1} is concluded.

\qed

\vskip.5cm

In  the proof of Theorem \ref{TEO:2} we shall use the following Lemma.
%
%
\medskip

\begin{lemma}

\label{lemma:Hold}

Let $\si>0$.
Then there is $c>0$ such that, for any $\delta$ and for any $t \in \delta \mathbb N$, $t\ge \si$,
\begin{equation}
\label{a4}
|S^{(\delta,+)}_t(u)(r)-S^{(\delta,+)}_t(u)(r')|\le c \max\{|r-r'|^\frac 1 3, \sqrt\delta\}
\end{equation}

\end{lemma}

\medskip

\noindent
{\bf Proof.}  It is clearly sufficient to bound the left hand side of \eqref{a4} when $|r-r'|$
and $\delta$ are such that:
\[
2\delta <\si, \qquad 2|r-r'|^{2/3} <\si
\]
We first consider the case when $|r-r'|^{2/3}\ge \delta$. We then have
\[
1\le \frac{|r-r'|^{2/3}}{\delta} \le \frac{\si}{\delta}-1
\]
Then there exists a positive integer $k^*$
such that $k^*\delta \le \si$ and
\[
\frac{|r-r'|^{2/3}}{\delta} \le k^* \le \frac{|r-r'|^{2/3}}{\delta} + 1
\]
We then apply \eqref{e6.3.0.1.2} with $s= t - k^*\delta$ getting
        \begin{equation}
        \label{r1}
| S^{(\delta,+)}_{t}(u)(r)- S^{(\delta,+)}_{t}(u)(r')| \le c'\Big( |r-r'|^{1/3}+   \sqrt{|r-r'|^{ 2/3}+\delta}\Big)\le c'(1+\sqrt 2) |r-r'|^{ 1/3}
         \end{equation}
Suppose next $|r-r'|^{2/3}\le \delta$. Choose $s=t-\delta$ then  \eqref{e6.3.0.1.2} gives:
\begin{equation}\label{r2}
| S^{(\delta,+)}_{t}(u)(r)- S^{(\delta,+)}_{t}(u)(r')| \le 2 c' \sqrt{\delta}
\end{equation}
so that \eqref{a4} follows from \eqref{r1} and \eqref{r2}. \qed

\vskip.5cm

We conclude the section with a corollary of the proof of Theorem \ref{thmee7.1}.

\vskip.5cm

\begin{theorem}
\label{thmee7.1.1}
Let $u\in C([0,1],\mathbb R_+)$, $F(0;u) > 0$.  Then
for any $\zeta>0$ there are $\tau_\zeta>0$
and $d_\zeta>0$ so that for any $\delta: \;F(0;u) > j\delta$,
for any $t\in \delta\mathbb N$, for any $t' \in \delta\mathbb N$, $t'\in (t,t+\tau_\zeta)$ and for any $r$ and $r'$ such that $|r-r'|<d_\zeta$,
 	\begin{equation}
	\label{e6.1aaa.1.11}
| S^{(\delta,+)}_{t}(u)(r)-S^{(\delta,+)}_{t}(u)(r')| <\zeta,
\quad| S^{(\delta,+)}_{t'}(u)(r)-S^{(\delta,+)}_{t}(u)(r)| <\zeta
	\end{equation}

\end{theorem}

\medskip

\noindent
{\bf Proof.}   It follows from \eqref{e4.3}
and the continuity of $u$ that for any $\zeta$ there is $d^*_\zeta$ so that
for any $t\ge 0$
 	\begin{equation}
	\label{e6.1aaa.1.12}
| G^{\rm neum}_t * u(r)- G^{\rm neum}_t * u(r')| <\frac{\zeta}2,
\quad|r-r'| <d^*_\zeta
	\end{equation}
Recalling \eqref{e6.3.0.1.2.1} we then set
 	\begin{equation}
	\label{e6.1aaa.1.13}
d_\zeta < \min\Big\{ d^*_\zeta; \;\;\;  \frac {\zeta^3}{4c'(2c')^2}\;;\;\frac{\zeta^3}{c''(2c')^2} \Big\}
	\end{equation}
As in the proof of Theorem \ref{thmee7.1} we first consider the case when $(2c')^2 \delta <\zeta^2$.  We then
choose $s<t$ as the smallest time in $\delta\mathbb N$ such that $2c'\sqrt{t-s} <\zeta$;   in the present case
where $t$ is not bounded away from 0 it may happen that $s=0$; if not the analysis is just as in the
proof of Theorem \ref{thmee7.1}. If instead $s=0$ we use \eqref{e6.1aaa.1.12}
to replace the bound in \eqref{e6.3} with $s=0$.  Then we can replace
\eqref{e6.3.0.1.2} by
 	\begin{equation}
	\label{e6.1aaa.1.14}
| S^{(\delta,+)}_{t}(u)(r)- S^{(\delta,+)}_{t}(u)(r')| \le
\frac{\zeta}2 +  c'\sqrt{t} < \zeta
	\end{equation}
The proof for the case when   $(2c')^2 \delta \ge \zeta^2$  is just as in the proof of Theorem \ref{thmee7.1} so that the first inequality in \eqref{e6.1aaa.1.11} is proved.

The second inequality in \eqref{e6.1aaa.1.11} follows from the first one by the same argument used
in the proof of Theorem \ref{thmee7.1} and since the first one has been proved without restrictions on $t$
the second one has also no restriction in $t$.  \qed

%

\vskip2cm

\section{Hydrodynamic limit}
\label{ottava}

{\bf Proof of Theorem \ref{Teo:EU}}. We
fix an element $u\in L^\infty([0,1],\mathbb R_+)$
such that $F(0;u)>0$.  We first restrict to
$\delta\in \Delta_\tau:= \{2^{-n}\tau, n \in  \mathbb N\}$, $\tau>0$
and prove convergence of $S^{(\delta,+)}_t(u)$
as $\delta\to 0$ in $\Delta_\tau$ when $t$ is
restricted to the interval $[\si,S]$,
$0<\si <S$.  More precisely we define a function $\psi^{(n)}(r,t)$ on $[0,1]\times [\si,S]$
by setting
\[
\psi^{(n)}(r,t) = S^{(2^{-n}\tau,+)}_t(u)(r),\quad r\in [0,1],\;\; t \in [\si,S] \cap (2^{-n}\tau)\mathbb N
\]
and defining $\psi^{(n)}(r,t)$ when $t \in [\si,S]$ by linear interpolation.

By Theorem \ref{thmee7.1} the family $\{\psi^{(n)}\}$ is equibounded
and equicontinuous hence by the Ascoli-Arzel\`a theorem it converges in
sup norm by subsequences to a continuous function $\psi(r,t)$  on $[0,1]\times [\si,S]$.
On the other hand for any $r\in [0,1]$ and $t\in [\si,S]\cap  \{k2^{-n}\tau,\; n,k\in \mathbb N\}$:
	\begin{equation*}
\lim_{m\to \infty} F(r; S^{(2^{-m}\tau,+)}_t(u)) = F(r;\psi(\cdot,t))
 	\end{equation*}
because, by \eqref{ee8.4}, $F(r; S^{(2^{-m}\tau,+)}_t(u))$ is a non increasing
function of $m$ which thus converges as $m\to \infty$. Thus all limit functions
$\psi(r,t)$ agree on $t\in [\si,S]\cap  \{k2^{-n}\tau,\; n,k\in \mathbb N\}$
and since they are continuous they agree on the whole $[\si,S]$,
thus the sequence $\psi^{(n)}(r,t)$ converges
in sup-norm  as $n\to \infty$ to a continuous function
$\psi(r,t)$.

By the arbitrariness of $\si$ and $T$ the function $\psi(r,t)$ extends to
the whole $[0,1]\times (0,\infty)$ and summarizing we have
	\begin{equation}
	\label{corr8.1}
\lim_{n\to \infty} \|S^{(2^{-n}\tau,+)}_t(u) - \psi(\cdot,t)\|_\infty  = 0,\quad
t>0, t\in  (2^{-n}\tau)\mathbb N
 	\end{equation}
the convergence
being uniform in $t\in \{ (2^{-n}\tau)\mathbb N\}$ when it varies
on the compacts not containing 0.

\vskip.5cm

\begin{proposition}
For any $r\in [0,1]$
	\begin{equation}
	\label{corr8.2}
\lim_{t\to 0} F(r;\psi(\cdot,t))  = F(r;u)
 	\end{equation}

\end{proposition}

\medskip

\noindent
{\bf Proof.}  Let $t= k 2^{-n}\tau$, $k$ and $n$ positive integers.  Then
by \eqref{ee8.4}
	\begin{equation*}
F(r;\psi(\cdot,t)) = \lim_{n\to \infty} F(r; S^{(2^{-n}\tau,+)}_t(u)) \le
F(r; S^{(t,+)}_t(u))
 	\end{equation*}
Let $X>0$, $r_X:= \max\{r-X,0\}$, then
	\begin{equation*}
F(r; S^{(t,+)}_t(u)) \le F(r_X,u) + F(0;u) \sup_{r'} \int_{|r-r'|>X} G^{\rm neum}_t(r,r')\, dr
 	\end{equation*}
By \eqref{e6.1aaa.01.1}
	\begin{equation*}
F(r; S^{(t,+)}_t(u)) \le F(r;u) + \|u\|_\infty X + F(0;u) \sqrt 2 e^{- X^2/(4t)}
	\end{equation*}
By choosing $X = t^{1/4}$
	\begin{equation*}
F(r; S^{(t,+)}_t(u)) \le F(r;u) + \|u\|_\infty \Big( t^{1/4}+  \sqrt 2 e^{- t^{-1/2}/4}\Big)
	\end{equation*}
To prove a lower bound we write
	\begin{equation*}
F(r;\psi(\cdot,t)) = \lim_{n\to \infty} F(r; S^{(2^{-n}\tau,+)}_t(u)) \ge
F(r; S^{(t,-)}_t(u))
 	\end{equation*}
and have
	\begin{equation*}
F(r; S^{(t,-)}_t(u)) \ge F(r+X,u) - F(0;u) \sup_{r'} \int_{|r-r'|>X} G^{\rm neum}_t(r,r')\, dr
 	\end{equation*}
	\begin{equation*}
F(r; S^{(t,-)}_t(u)) \ge F(r;u) - \|u\|_\infty \Big( t^{1/4}+  \sqrt 2 e^{- t^{-1/2}/4}\Big)
	\end{equation*}
Thus
	\begin{equation*}
 |F(r;\psi(\cdot,t))  - F(r;u)| \le \|u\|_\infty \Big( t^{1/4}+  \sqrt 2 e^{- t^{-1/2}/4}\Big),
 \quad t= k 2^{-n}\tau>0
	\end{equation*}
By the continuity of $\psi(\cdot,t)$ and because the set $\{ k 2^{-n}\tau$, $k\in \mathbb N_+, n\in \mathbb N\}$
is dense in $\mathbb R_+$, it follows that
	\begin{equation*}
\sup_{t \le S} |F(r;\psi(\cdot,t))  - F(r;u)| \le \|u\|_\infty \Big( S^{1/4}+  \sqrt 2 e^{- S^{-1/2}/4}\Big)
	\end{equation*}
hence \eqref{corr8.2}.  \qed

\vskip.5cm

\begin{proposition}
For any $t\in \{ k 2^{-n}\tau$, $k\in \mathbb N_+, n\in \mathbb N\}$,
 	\begin{equation}
	\label{corr8.3}
\lim_{n\to \infty} \int dr|\psi(r,t)-S^{(2^{-n}\tau,-)}_t(u)(r)| = 0
 	\end{equation}
  	\begin{equation}
	\label{corr8.4}
F(r;\psi(\cdot,t))\ge F(r;S^{(2^{-n}\tau,-)}_t(u)),\quad r\in [0,1]
 	\end{equation}

\end{proposition}

\medskip

\noindent
{\bf Proof.}
\eqref{corr8.3} follows from \eqref{corr8.1} and \eqref{e6.1aaa.2}. By
\eqref{corr8.3}
  	\begin{equation*}
F(r;\psi(\cdot,t))=\lim_{n\to \infty} F(r;S^{(2^{-n}\tau,-)}_t(u))
 	\end{equation*}
which implies \eqref{corr8.4} because, by \eqref{ee8.4}, $F(r;S^{(2^{-n}\tau,-)}_t(u))$ is a non decreasing function of $n$.  \qed

\vskip.5cm
By \eqref{ee8.4} we then have for all $r\in [0,1]$ and all $\delta$ and $t$ in
$\{ k 2^{-n}\tau$, $k\in \mathbb N_+, n\in \mathbb N\}$,
  	\begin{equation}
	\label{corr8.5}
F(r;\psi(\cdot,t))\ge F(r;S^{(\delta,-)}_t(u)),\quad
F(r;\psi(\cdot,t))\le F(r;S^{(\delta,+)}_t(u))
 	\end{equation}
\eqref{corr8.5} does not yet prove that $\psi$ separates the barriers
because we have to consider all $t$ and $\delta$ and not only those above.
To this end we observe that the function $\psi(r,t)$ that we have defined so far actually depends
on the initial choice of $\tau$, to make this explicit we write $\psi_\tau(r,t)$.  Of course
we have for all $\tau>0$:
  	\begin{equation}
	\label{corr8.6}
F(r;S^{(\delta,-)}_t(u)) \le
F(r;\psi_{\tau}(\cdot,t)) \le F(r;S^{(\delta,+)}_t(u)),
\quad \delta, t\in
\{ k 2^{-n}\tau, k\in \mathbb N_+, n\in \mathbb N\}
 	\end{equation}
so that we only need to show that $\psi_\tau$ does not depend on $\tau$.
To prove independence of $\tau$ we use the following lemma:

%
%
%
%

\vskip.5cm

   \begin{lemma}
  \label{lemmae35}
There is $c$ so that for any  $0<\delta<\delta'$,
$u\in \mathcal U_\delta$  and $n\ge 1$
  \begin{equation}
    \label{e6.20.0}
| S^{(\delta,-)}_{n\delta}(u) -S^{(\delta',-)}_{n\delta'}(u) |_1 \le c |u|_1 n\frac{\delta'-\delta}{\delta^{3/2} }
   \end{equation}

    \end{lemma}

\medskip

\noindent
{\bf Proof.} In order to compare $S^{(\delta,-)}_{\delta}$ and
$S^{(\delta',-)}_{\delta'}$ we shall use the following bounds:
  \begin{equation}
    \label{e6.20.1}
| K^{(\delta)}(w) -K^{(\delta')}(w) |_1 \le 2j(\delta'-\delta),\quad
|G^{\rm neum}_\delta*w - G^{\rm neum}_{\delta'}*w|_1 \le \frac{c(\delta'-\delta)}{\delta^{3/2} } |w|_1
   \end{equation}
together with $| K^{(\delta)}(w) -K^{(\delta)}(w) |_1 \le |v-w|_1$, see \eqref{e6.17}.
Indeed we can bound $| S^{(\delta,-)}_{\delta}(w) -S^{(\delta',-)}_{\delta'}(v) |_1$ by
  \begin{eqnarray*}
&&\le |K^{(\delta)}\{G^{\rm neum}_\delta*w-G^{\rm neum}_{\delta'}*v\}|_1 + |(K^{(\delta')}-K^{(\delta)})G^{\rm neum}_{\delta'}*v\}|_1 \\
&&  \le |G^{\rm neum}_\delta*w-G^{\rm neum}_{\delta'}*v|_1 + 2j(\delta'-\delta)
\\
&& \le |G^{\rm neum}_\delta*w-G^{\rm neum}_{\delta}*v|_1 +|G^{\rm neum}_\delta*v-G^{\rm neum}_{\delta'}*v|_1 + 2j(\delta'-\delta)
   \end{eqnarray*}
getting
 \begin{equation}
     \label{e6.20.1.1}
     | S^{(\delta,-)}_{\delta}(w) -S^{(\delta',-)}_{\delta'}(v) |_1
\le |w-v|_1
+ c\frac{\delta'-\delta}{\delta^{3/2} } |v|_1  + 2j (\delta'-\delta)
   \end{equation}
By using \eqref{e6.20.1.1} with $w=S^{(\delta,-)}_{(n-1)\delta}(u)$ and
$v=S^{(\delta',-)}_{(n-1)\delta'}(u)$,  then, by iteration, we get \eqref{e6.20.0}.  \qed

\vskip.5cm

   \begin{theorem}
  \label{thme35}
$\psi_ \tau$ is independent of $\tau$.

    \end{theorem}

\medskip

\noindent
{\bf Proof.}  We shall prove that for any $\tau$ and $\tau'$
\[
F(r;\psi_\tau(\cdot,t)) = F(r;\psi_{\tau'}(\cdot, t)),\quad r\in [0,1], \;\;t > 0
\]
and this will prove Theorem \ref{thme35}.  We suppose that
$\tau' \notin \{ k\tau 2^{-n}, k, n \in \mathbb N\}$ (otherwise the
statement trivially holds). We fix $t'= n \delta'$, $\delta'= \tau' 2^{-m}$.
Let $\delta = k \tau 2^{-q}$, $\delta <\delta'$.
By the previous lemma, for all $r\in [0,1]$
 \begin{equation*}
F(r; S^{(\delta',-)}_{t'}(u) )\le F(r;S^{(\delta,-)}_{n\delta}(u)) +
c F(0;u) n\frac{\delta'-\delta}{\delta^{3/2} }
   \end{equation*}
Write $\delta= k_p \tau 2^{-p}$ so that $k_p = k 2^{p-q}$ is a positive integer for $p$ large enough.
Then by \eqref{ee8.4}
 \begin{equation*}
 F(r;S^{(\delta,-)}_{n\delta}(u) )\le  F(r;S^{(\tau 2^{-p},-)}_{n\delta}(u))
   \end{equation*}
By taking $p\to \infty$:
 \begin{equation*}
F(r; S^{(\delta',-)}_{t'}(u) )\le F(r;\psi_{\tau}(\cdot,n\delta)) +
c F(0;u) n\frac{\delta'-\delta}{\delta^{3/2} }
   \end{equation*}
We then let $\delta \to \delta'$ on $\{k\tau 2^{-n}, k, n \in \mathbb N\}$.  In this limit
$n\delta \to t'$ and by the continuity of $\psi_\tau(\cdot, s)$ in $s$ we get
 \begin{equation*}
F(r; S^{(\delta',-)}_{t'}(u) )\le F(r;\psi_\tau(\cdot, t'))
   \end{equation*}
We next take $m\to \infty$, recall  $\delta'= \tau' 2^{-m}$, and get
 \begin{equation*}
F(r; \psi_{\tau'}(\cdot,t')) \le F(r;\psi_\tau(\cdot, t')),\qquad \text{for any $t'\in \{k\tau' 2^{-n}, k, n \in \mathbb N\}$}
   \end{equation*}
In an analogous fashion we get
 \begin{equation*}
F(r;\psi_\tau(\cdot, t))\le F(r; \psi_{\tau'}(\cdot,t)) ,\qquad \text{for any $t\in \{k\tau 2^{-n}, k, n \in \mathbb N\}$}
   \end{equation*}
Then $\psi_\tau(\cdot,t)= \psi_{\tau'}(\cdot,t)$ for all $t$ in a dense set, hence they are equal everywhere
being both continuous.\qed

\bigskip

\noindent
The proof of Theorem \ref{Teo:EU} is concluded.
\qed

\bigskip

\noindent
{\bf Proof of Theorem \ref{Teo:Char}}. It follows from the reasoning above and the use of Theorem \ref{thmee8.1} with the choice $m=2^n$
and $\delta=t2^{-n}$.
\qed

\vskip.5cm
\noindent
{\bf Proof of Theorem \ref{Teo:Hydro}.}
 The proof of Theorem \ref{Teo:Hydro} is an immediate consequence of Theorem \ref{thme35}.
\qed

\vskip.5cm

\medskip
\noindent
We are left with  the proof of Theorem \ref{TEO:2}, that we explain in the remaining part of this section.
We fix $\rho_{\rm init}$ such that $\rho_{\rm init}(1)>0$ and we call $\rho_t$ the function of
Theorem \ref{Teo:Hydro}.

For any $a>0$ arbitrarly small we define
$$
T_a = \sup\{t >0 :  \rho_t(1) \ge a\}
$$

		\begin{lemma}
\label{lemma:G1}

For any $a>0$ there exists $0<a'<a$ such that

		\begin{equation}
		\label{G2}
S^{(\delta,+)}_{n\delta}
\rho_{\rm init}(1)\ge a'
\qquad \text{for any} \quad n \quad \text{such that} \quad \delta n < T_a
	\end{equation}
	\end{lemma}
\medskip

\noindent

{\bf Proof.} Let $t \in \delta \mathbb N$ with $t < T_a$ then $\rho_t(1)\ge a$.
From Theorem \ref{Teo:Char}, for any $r \in [0,1]$,
$t \in \delta \mathbb N$ we have
\begin{equation}
F(r; S_t^{(\delta,-)}(\rho_{ \rm init})) \le F(r;\rho_t) \le F(r; S_t^{(\delta, +)}(\rho_{\rm init}))
\end{equation}
On the other hand, from \eqref{e6.1aaa.2}, for any $r \in [0,1]$, $t \ge 0$,
\begin{equation}
\big|F(r; S_t^{(\delta,-)}(\rho_{\rm init}))-F(r; S_t^{(\delta,+)}(\rho_{ \rm init}))\big|\le 4 j \delta\;.
\end{equation}
 As a consequence, writing $\rho^{(\delta,+)}_t := S_t^{(\delta,+)}(\rho_{ \rm init})$ and choosing $r=1-\sqrt{\delta}$, we have
    \begin{equation}
    \label{G1}
\int_{1-\sqrt \delta}^1 \rho^{(\delta,+)}_t(r) \, dr\ge \int_{1-\sqrt \delta}^1 \rho_t(r) \, dr- 4 j \delta
\end{equation}
From Lemma \ref{lemma:Hold}, for $r \in [1-\sqrt \delta, 1]$,
	\begin{equation*}
|\rho^{(\delta,+)}_t(r)-\rho^{(\delta,+)}_t(1)|\le c \max\{|1-r|^{\frac 1 3}, \sqrt{\delta}\} \le c \delta^{\frac 1 6}
	\end{equation*}
hence
\begin{equation}
\label{8888}
\int_{1-\sqrt \delta}^1 \rho^{(\delta,+)}_t(r) \, dr\le (\rho^{(\delta,+)}_t(1) + c \, \delta^{\frac 1 6})\,  \sqrt \delta
\end{equation}
Combining \eqref{G1} and \eqref{8888} we have
\begin{equation}
\rho^{(\delta,+)}_t(1) + c \, \delta^{\frac 1 6} \ge \frac 1 {\sqrt \delta} \int_{1-\sqrt \delta}^1 \rho^{(\delta,+)}_t(r) \, dr \ge \frac 1 {\sqrt \delta}  \int_{1-\sqrt \delta}^1 \rho_t(r) \, dr- 4j \, \sqrt\delta
\end{equation}
thus
\begin{equation}
\label{9999}
\rho^{(\delta,+)}_t(1) \ge \frac 1 {\sqrt \delta}  \int_{1-\sqrt \delta}^1 \rho_t(r) \, dr- c' \, \delta^{\frac 1 6}
\end{equation}
From the space continuity of $\rho_t$ obtained in Theorem \ref{Teo:Hydro}, for any $a>0$ there exists $\delta>0$ small enough such that, for $|r-1|\le \sqrt \delta$,
$$\rho_t(r)\ge \rho_t(1)- a/2\ge a/2
$$
where the last inequality holds for all $t< T_a$. Then the statement of the Lemma follows
from \eqref{9999} with  $a'=a/2-c' \delta^{\frac 1 6}$ which is positive for $\delta$ small enough.
\qed

\vskip1cm

\begin{lemma}\label{lemma:G2}
For any $a>0$ there is $C_a>0$ such that for any  $t  \in \delta  \mathbb N$, $t  < T_a$
\begin{equation}\label{G3}
R_\delta (S^{(\delta ,+)} _{n\delta}{(\rho_{\rm init})})\ge 1-C_a\, \delta
\end{equation}
\end{lemma}
\medskip
\noindent
{\bf Proof.} Fix $C>0$  and denote $\rho^{(\delta,+)}_t := S_t^{(\delta,+)}(\rho_{\rm init})$. From Lemma \ref{lemma:Hold} we know that there is $c>0$ so that for any $r \in [1-C \delta, 1]$, $t \in \delta \mathbb N$,
\begin{equation}
\rho^{(\delta,+)}_t(r)\ge  \rho^{(\delta,+)}_t(1)- c  \delta^\frac 1 3
\end{equation}
then, from Lemma \ref{lemma:G1}, for any $a>0$ there is $0<a'<a$ such that
\begin{equation}
\int_{1-C\delta}^1 \rho^{(\delta,+)}_t(r) \, dr \ge C \delta (a' -c  \delta^{\frac 1 3}) \qquad \forall t< T_a
\end{equation}
now it is sufficient to chose $C=C_a> (a' -c  \delta^{\frac 1 3})/j$, $\delta$ small enough to get
\begin{equation}
\int_{1-C_a\delta}^1 \rho^{(\delta,+)}_t(r) \, dr > j\delta \qquad \forall t< T_a,
\end{equation}
that gives \eqref{G3}. \qed

\vskip1cm

\noindent
{\bf Proof of Theorem \ref{TEO:2}}.
We define the dynamics

\begin{eqnarray}\label{G40}
\hat S_{n\delta}^{(\delta,+)}(u) &:=& G_{\delta}^{\rm neum} * \dots * Q^\delta G_{\delta}^{\rm neum} *  Q^\delta u \qquad n \text{ times} \nonumber\\
&=& G_{\delta}^{\rm neum} * Q^\delta \hat S_{(n-1)\delta}^{(\delta,+)}(u)
\end{eqnarray}
with
\begin{equation}
Q^\delta u = u +j \delta D_0 -j \delta D_1
\end{equation}
then
\begin{equation}
\hat S_{n\delta}^{(\delta,+)}(u)= G_{n \delta}^{\rm neum} * u + j \delta \sum_{k=0}^{n-1}G_{k \delta}* D_0 - j \delta \sum_{k=0}^{n-1}G_{k \delta}* D_1
\end{equation}
hence $\hat S_t ^{(\delta,+)}(u)$ converges as $\delta \to 0$ to the dynamics defined by \eqref{e2.5}.
It remains to prove that $S_{n\delta}^{(\delta,+)}(u)- \hat S_{n\delta}^{(\delta,+)}(u)$ converges weakly to zero for $n \delta <\sup_a T_a$.

From \eqref{G40} and \eqref{e2.10}, we can write
\begin{eqnarray}
&&  S_{n\delta}^{(\delta,+)}(u)- \hat S_{n\delta}^{(\delta,+)}(u)= \nonumber\\
&& = G_{\delta}^{\rm neum} *\left(K^\delta  S_{(n-1)\delta}^{(\delta,+)}-Q^\delta \hat S_{(n-1)\delta}^{(\delta,+)}\right)(u)= \nonumber\\
&&  = G_{\delta}^{\rm neum} *\left(K^\delta-Q^\delta \right) S_{(n-1)\delta}^{(\delta,+)}(u) +  G_{\delta}^{\rm neum} * Q^\delta \left(S_{(n-1)\delta}^{(\delta,+)}-\hat S_{(n-1)\delta}^{(\delta,+)}\right)(u)= \nonumber\\
&&  = G_{\delta}^{\rm neum} *\left(K^\delta-Q^\delta \right) S_{(n-1)\delta}^{(\delta,+)}(u) +G_{\delta}^{\rm neum} * Q^\delta G_{\delta}^{\rm neum}* \left(K^\delta S_{(n-2)\delta}^{(\delta,+)}-Q^\delta \hat S_{(n-2)\delta}^{(\delta,+)}\right)(u)= \nonumber\\
&&= \sum_{k=1}^n G_{\delta}^{\rm neum} *  Q^\delta G_{\delta}^{\rm neum} * \dots * Q^\delta G_{\delta}^{\rm neum} *(K^\delta-Q^\delta)S_{(n-k)\delta}^{(\delta,+)} (u) \qquad (\text{by iteration}) \label{777}
\end{eqnarray}
where the $G_{\delta}^{\rm neum}$ appears $k$ times in the $k$-th term of the sum and
\begin{equation}
(K^\delta -Q^\delta) v :=    j \, \delta \, D_1- \mathbf 1_{[R_{\delta}(v),1]} \, v
\end{equation}
Then, in order to prove the convergence of \eqref{777} to 0  we  prove that each term in the sum \eqref{777} converges to 0 as $\delta \to 0$. This is true since   for any  $n  : \; n \delta < \sup_a T_a$
\begin{equation}
(K^\delta- Q^\delta)S^{(\delta,+)}_{n \delta}u \to 0 \qquad \text{weakly as } \delta \to 0
\end{equation}
The proof of this last statement follows form the following argument.
We first fix $a>0$ arbitrarily small, then, from   \eqref{G3}, there exists $C_a>0$ so that
\begin{equation}
\bigg|\text{supp} \left(\mathbf 1_{R_{\delta}(S^{(\delta,+)}_{k \delta}u)} \, S^{(\delta,+)}_{n \delta}u \right) \bigg| \le C_a \delta, \qquad \text{for any } \; n : \: n\delta \le T_a
\end{equation}
Then
for any test function $\phi$, $n\delta \le T_a$,
\begin{eqnarray}
&& \Big| \frac 1 {j \delta} \int_{R_{\delta}(S^{(\delta,+)}_{n \delta}u)}^1 \, S^{(\delta,+)}_{n\delta}u(r) \, \cdot \phi(r) \, dr - \phi(1) \Big| \nonumber \\
&& =\Big| \frac 1 {j \delta} \int_{R_{\delta}(S^{(\delta,+)}_{n \delta}u)}^1 \, S^{(\delta,+)}_{n\delta}u(r) \, \cdot \left(\phi(r)-\phi(1)\right) \, dr  \Big| \nonumber
\\
&& \le \sup_{r \in [R_{\delta}(S^{(\delta,+)}_{n \delta}u),1]} \big| \phi(r)-\phi(1) \big| \le  \sup_{|r-1| \le C_a \delta} \big| \phi(r)-\phi(1) \big|
\end{eqnarray}
that vanishes as $\phi$ is  continuous. Hence, for any $a>0$,
\begin{equation}\label{G4}
\lim_{\delta \to 0} \Big| \frac 1 {j \delta} \int_{R_{\delta}(S^{(\delta,+)}_{n \delta}u)}^1 \, S^{(\delta,+)}_{n\delta}u(r) \, \cdot \phi(r) \, dr - \phi(1) \Big| =0 \qquad \text{for }\: k \delta \le T_a
\end{equation}
 then \eqref{G4} is certainly true as long as $n \delta \le \sup_a T_a$, this yields the convergence in  distribution to equation \eqref{e2.5}  for any time $t$ such that $\rho_t(1)>0$.  We know that the convergence of $S_t^{(t 2^{-n},+)}(\rho_{\rm init})$ to $\rho_t$ as $n \to \infty$ in the sense of the interfaces (see Theorem \ref{Teo:Char}) implies weak convergence against smooth test functions. This and the uniqueness of the weak limit  univocally characterizes $\rho_t$ as the function given by \eqref{e2.5} for $t$  such that $\rho_t(1)>0$. Then the Theorem is proved.
\qed

\newpage
{\bf Acknowledgments.}
We thank Pablo Ferrari and John Ockendon for many useful comments and discussions.
The research has been partially supported by PRIN 2009 (prot. 2009TA2595-002)
 and FIRB 2010 (grant n. RBFR10N90W).
A. De Masi and E. Presutti acknowledge kind hospitality at the Dipartimento di Matematica
della Universit\`a di Modena.
 G. Carinci and C. Giardin\`a thank Universit\`a dell'Aquila for  welcoming during their
visit at Dipartimento di Matematica.

\end{document}